\documentclass[a4paper,11pt, reqno]{amsart}
\usepackage[DIV=12, oneside]{typearea} 
\usepackage{ucs}
\usepackage[utf8x]{inputenc}
\usepackage[T1]{fontenc}
\usepackage[english]{babel}
\usepackage{esint}
\usepackage{amssymb, amsthm, bbm}
\usepackage{enumerate, mathrsfs, hyphenat, graphicx}
\usepackage[font=sf, labelfont=sf]{caption}
\usepackage{thmtools, thm-restate}
\usepackage{thm-patch, thm-kv, thm-autoref}
\usepackage[colorlinks=true, urlcolor=blue, linkcolor=blue, citecolor=black]{hyperref}
\theoremstyle{plain}

\declaretheorem[title=Theorem]{theorem}
\declaretheorem[title=Lemma,numbered=no]{lemma*}
\declaretheorem[title=Lemma,sibling=theorem]{lemma}
\declaretheorem[title=Proposition,sibling=theorem]{proposition}
\declaretheorem[title=Definition,sibling=theorem]{definition}
\declaretheorem[title=Corollary,sibling=theorem]{corollary}

\usepackage[backgroundcolor=white, bordercolor=blue,
linecolor=blue]{todonotes}

\newcommand{\R}{  \mathbb{R}}

\newcommand{\N}{  \mathbb{N}}

\renewcommand{\P}{ \mathbb{P}}

\newcommand{\E}{  \mathbb{E}}

\newcommand{\cB}{  \mathcal{B}}

\newcommand{\cS}{  \mathcal{S}}

\renewcommand{\epsilon}{\varepsilon}

\renewcommand{\d}{\, \textnormal{d}}

\makeatletter 
\def\Links{\tagsleft@true}\def\Rechts{\tagsleft@false} 
\makeatother

\usepackage{array}
\usepackage{booktabs}
\usepackage{stmaryrd}

\newcommand{\opA}{A}

\begin{document}
\renewcommand{\sectionautorefname}{Section}

\title{Intrinsic scaling properties for nonlocal operators}
\author{Moritz Kassmann, Ante Mimica}
\begin{abstract}
We study integrodifferential operators and regularity estimates 
for solutions to integrodifferential equations. Our 
emphasis is on kernels with a critically low singularity which does not allow for standard 
scaling. For example, we treat operators that have a logarithmic order of 
differentiability. For corresponding equations we prove a growth lemma and derive a priori 
estimates. We derive these estimates by classical methods developed for partial 
differential 
operators. Since the integrodifferential operators under consideration generate Markov 
jump processes, we are able to offer an alternative approach using probabilistic 
techniques.  
\end{abstract}

\address{Fakult\"{a}t f\"{u}r Mathematik\\
Universit\"{a}t Bielefeld \\
Postfach 100131 \\
D-33501 Bielefeld}
\email{moritz.kassmann@uni-bielefeld.de}

\address{Department of Mathematics\\
University of Zagreb\\
Bijeni\v{c}ka c. 30\\
10000 Zagreb, Croatia}
\email{amimica@math.hr}

\thanks{Research supported by German Science Foundation (DFG) via SFB
701. Supported in part by Croatian Science Foundation under the project 
3526.}



\date{March 21, 2015}

\maketitle


\section{Introduction}\label{sec:intro-main}

In recent years, regularity results for linear and nonlinear integrodifferential 
operators have been addressed by many research articles. Scaling properties are crucially 
used in these approaches. We reconsider these cases and, at the same time, include limit 
cases where standard scaling properties do not hold anymore. We study
linear operators of the form 
\begin{align*}
 A u(x)&=\int\limits_{\R^d}\big(u(x+h)-u(x)-\langle
\nabla
u(x),h\rangle \mathbbm{1}_{B_1}(h)\big)K(x,h) \d h \,,
\end{align*}
which, provided certain assumptions on $K(x,h)$ are satisfied, are well defined 
for smooth and bounded functions $u:\R^d \to \R$. The quantity $K(x,h)$ equals 
the jump intensity of jumps from $x \in \R^d$ to 
$x+h \in \R^d$ for the Markov process $X$ that is generated by the linear 
operator $A$. If $K$ is independent of the first variable, then $X$ is a 
L\'{e}vy process. If $K(x,h) = |h|^{-\alpha -d} 
m(\frac{h}{|h|})$ for all $h \ne 0$ and some appropriate function 
$m:\mathbb{S}^{d-1} \to [0,+\infty]$, then the increments of $X$ have 
stable 
distributions. Looking at the 
operator $A$ as an integrodifferential operator, this property is important 
because it allows to use scaling techniques. 

Scaling techniques themselves are 
crucial when studying regularity properties of functions $u:\R^d \to \R$ 
satisfying the equation $Au=f$ in a domain $\Omega \subset \R^d$ for some 
function $f:\Omega \to \R$. In this work we study such properties with the 
emphasis on two features. We do not assume any regularity of the kernel 
function $K$ with respect to the first variable except for 
boundedness. Moreover, and this is the main new contribution, we systematically 
study classes 
of kernels that do not possess the aforementioned scaling property. 

Our results include a growth lemma (expansion of 
positivity) and Hölder-type regularity estimates. Moreover, we provide several estimates 
on the corresponding Markov jump process. Recall that, in the case of an 
elliptic operator of second order $Au = a_{ij}(\cdot) 
\partial_{i} \partial_{j} u$, the standard growth lemma reads as follows:

\begin{lemma}\label{lem:growth} There is a constant $\theta \in (0,1)$ such that,
if $R>0$ and $u:\R^d \to \R$ satisfying
\begin{align*}
 &-A u \leq 0 \text{ in } B_{2R} \,, \qquad u \leq 1 \text{ in }
B_{2R}, \qquad |(B_{2R} \!\setminus\! B_R) \cap \{u \leq 0\}| \geq \tfrac12
|B_{2R} \!\setminus\! B_R|\,,
\end{align*}
then $u \leq 1-\theta$ in $B_R$.
\end{lemma}

The above lemma holds true also for several nonlinear
operators. Such lemmas are systematically studied and applied in \cite{Lan71}.
Their importance is underlined in the article \cite{KrSa79}, in
which the authors establish a priori bounds for elliptic
equations of second order with bounded measurable coefficients. Nowadays they
form a standard
tool for the study of various questions of nonlinear partial differential
equations of second order, see \cite{CaCa95} and \cite{DBGV12}. Note that the
property
formulated in \autoref{lem:growth} is also referred to as expansion
of positivity which describes the corresponding property for $1-u$. 

In this work we prove a similar growth lemma for integrodifferential operators, see 
\autoref{lem:expans1} below. An important instance of an 
operator $A$ that we have in mind is
\begin{align}\label{eq:exa_L}
 A u(x)&=\int\limits_{\R^d} \big[ u(x\!+\!h) - u(x) \big]  a(x,h) |h|^{-d} 
\mathbbm{1}_{B_1}(h) \d h \,,
\end{align}
for some measurable function $a:\R^d \times \R^d \to [1,2]$. Note that, in the 
last years, similar results have been studied for kernels of the form $a(x,h) 
\asymp |h|^{-\alpha}$  
for some $\alpha \in (0,2)$ and we refer the 
reader to the 
short discussion 
below. The case $\alpha = 0$ is of 
particular interest because in 
this case the corresponding growth lemma fails. Our results apply to more 
general kernels than the one appearing in \eqref{eq:exa_L}. 

\subsection{Main assumptions and results} Let $K\colon 
\R^d\times(\R^d\setminus\{0\})\rightarrow 
[0,+\infty)$ be a measurable function such that $K(x,h)=K(x,-h)$ for all $x,h$  
and
\begin{align}
& \kappa^{-1}|h|^{-d}\ell(|h|)\leq K(x,h)\leq \kappa 
|h|^{-d}\ell(|h|)\quad \text{ for }\quad 0<|h|<R_0 \,, \label{eq:A1}  
\tag{A$_1$} \\
& \sup\limits_{x\in\R^d}\int\limits_{\R^d\setminus \{0\}}(1\wedge|h|^2)K(x,h) 
\d 
h \leq K_0 \,, \tag{K$_0$} \label{eq:K-Levy-bound}
\end{align}
where $\kappa, K_0 \geq 1$, $R_0\in (0,+\infty]$ are fixed constants and  
$\ell\colon 
(0,R_0)\rightarrow 
(0,+\infty)$ is a function satisfying $\int_0^{R_0} \ell(s)\frac{ds}{s} = + 
\infty$ and, for some $c_L \in (0,1)$, $c_U \geq 1$, 
and $\gamma\in (0,2)$, the following:
\begin{alignat}{2}
\label{eq:ell1} \tag{$\ell_1$} &\int_r^{R_0} \ell(s)\frac{ds}{s} < 
+ \infty \quad 
&& \text{ 
for } 0<r<R_0  \,, \\
\label{eq:ell2}  \tag{$\ell_2$} &\frac{\ell(r \lambda)}{\ell(r)} \geq c_L 
\lambda^{-\gamma}  
&& \text{ for } r>0 \text{ and } 1 \leq \lambda < \frac{R_0}{r}\,, \\
\label{eq:ell3}  \tag{$\ell_3$} &\frac{\ell(r \lambda )}{\ell(r)} \leq c_U 
\lambda^{d} && \text{ for } r>0 \text{ and }  1 \leq \lambda < 
\frac{R_0}{r} \,.
\end{alignat}
The last two conditions often are referred to as weak scaling conditions.

{\bf Examples:} The standard example is given by $\ell(s)=s^{-\alpha}$ for some 
$\alpha \in (0,2)$. Other examples satisfying the above 
conditions include $\ell(s) = 
s^{-\alpha} 
g(s)$ for $\alpha \in (0,2)$ and $g$ a function that varies slowly at 
$0$. More 
generally, the conditions \eqref{eq:ell2} and 
\eqref{eq:ell3} are satisfied if the function 
$\ell$ is regularly varying at zero of order $-\alpha \in (-2,0]$ and satisfies 
some weak bounds for large values of $s$.  The case $\alpha=0$ is 
very interesting. The choice $\ell(s) = 1$ is possible if $R_0 < + \infty$. In 
the case $R_0 = + \infty$ an interesting example is provided by $\ell(s) = 
\mathbbm{1}_{(0,1)}(s)+s^{-\gamma}\mathbbm{1}_{[1,+\infty)}(s)$ for some 
$\gamma > 0$.

We define an auxiliary function $L\colon (0,R_0)\rightarrow (0,+\infty)$ by 
$L(r)=\int_r^{R_0}\frac{\ell(s)}{s}\d s$, which is strictly decreasing. 
Note that, under our assumptions $L(0+)=+\infty$ and $L(R_0-)=0$. Furthermore, 
we define a measure $\mu$ on $\cB(B_{R_0})$ by
\begin{align*}
 \mu(dy)=\frac{\ell(|y|)}{L(|y|)}\frac{dy}{|y|^d} 
\end{align*}
and, for $a>1$, a scale function $\varphi_a=\varphi\colon (0,R_0)\rightarrow 
(0,\infty)$ by 
$\varphi(r)=L^{-1}(a^{-1}L(r))$. 

Define an operator $A\colon C^2_b(\R^d)\rightarrow C(\R^d)$ by
\begin{align}\label{eq:def_A}
 A u(x)&=\int\limits_{\R^d\setminus \{0\}}\big(u(x+h)-u(x)-\langle
\nabla
u(x),h\rangle \mathbbm{1}_{B_1}(h)\big)K(x,h) \d h 
\end{align}
where $K\colon \R^d\times \R^d\setminus\{0\}\rightarrow [0,+\infty)$ satisfies 
\eqref{eq:A1} and \eqref{eq:K-Levy-bound}.

Now we can formulate our first main result, a growth lemma for nonlocal 
operators. We state the result for functions which, together with their first 
and second 
derivatives, are continuous and bounded. It is an important feature of this result, that 
none of the constants depends on the regularity of the function under consideration. 
Thus, the result is tailored for later 
applications to viscosity solutions of fully nonlinear partial differential equations. 

\begin{lemma}\label{lem:expans1}  Assume \eqref{eq:K-Levy-bound} and 
\eqref{eq:A1} hold true with $R_0 = +\infty$. Let $\eta, \delta \in 
(0,1)$ and $C_0>0$.
There exist constants $a>2$ and $\theta\in (0,1)$ such that, if for $r > 0$ and 
$v \in C^2_b(\R^d)$ 
the following conditions are satisfied:
\begin{alignat*}{2}
 - Av (x) &\leq L(\varphi(r)) \quad
  (=a^{-1}L(r)) \qquad  &&\text{ for } x \in
B_{\varphi(r)} \,,
\\
 v (x) &\leq 1 &&\text{ for } x\in B_{\varphi(r)} \,, \\
 v (x) &\leq \, C_0 \left(\frac{L(\varphi(r))}{L(|x|)}\right)^\eta &&\text{ for 
} x 
\in\R^d\setminus
B_{\varphi(r)} \,, \\
 \mu\big((B_{\varphi(r)}\setminus B_r) \cap \{v \leq 0\}\big)  &\geq \delta 
\mu\big(B_{\varphi(r)}\setminus B_r \big)\,,
\end{alignat*}
then the following is true:
\begin{equation}\label{eq:assertion_lemma}
 v(x) \leq 1-\theta \quad \text{ for all }\quad  x \in B_{r}. 
\end{equation}
\end{lemma}

{\bf Remark:} As the proof shows, the value of $\theta$ is a multiple of 
$a^{-1}$. 

Note that, in the by now well-known case where $\ell(r) = r^{-\alpha}$ and $L(r) \asymp 
r^{-\alpha}$ for $\alpha \in (0,2)$, this result reduces to a growth lemma which is 
very similar to those given in \cite{Sil06} and \cite{CaSi09}. Let us now formulate the 
second main result.

\begin{theorem}\label{theo:main-analysis} Assume \eqref{eq:K-Levy-bound} and 
\eqref{eq:A1} hold true with $R_0 \in (0, +\infty]$. There exist 
constants $c>0$ and 
$\beta\in (0,1)$
such that for $0 < r \leq \tfrac{R_0}{2}$, $f \in L^\infty(B_r)$, and $u 
\in 
C^2_b(\R^d)$ satisfying $Au=f$ in 
$B_r$, the following holds
 \begin{align}\label{eq:theo-estim}
  \sup_{x,y\in B_{r/4}}\tfrac{|u(x)-u(y)|}{L(|x-y|)^{-\beta}}\leq 
c L(r)^{\beta}\|u\|_\infty + c L(r)^{\beta-1} \|f\|_{L^\infty(B_r)}\,.
 \end{align}
If $R_0 = + \infty$, then \eqref{eq:theo-estim} holds true for every $r > 
0$. 
\end{theorem}
\enlargethispage{3ex}

In the case $\ell(r) = r^{-\alpha}$, we set $\nu=\alpha \beta$, and the 
estimate 
\eqref{eq:theo-estim} reduces to
 \begin{align*}
  \sup_{x,y\in B_{r/4}}\tfrac{|u(x)-u(y)|}{|x-y|^{\nu}}\leq 
c r^{-\nu}\|u\|_\infty + c r^{\alpha - \nu} \|f\|_{L^\infty(B_r)}\,,
 \end{align*}
which one would expect from standard scaling behavior of the 
integrodifferential operator.

In the case $R_0=+\infty$ we obtain a Liouville theorem.
\begin{corollary}
If \eqref{eq:A1} holds for $R_0=+\infty$, then every function $u \in 
C^2_b(\R^d)$ satisfying 
$Au=0$ on $\R^d$ is a constant function. 
\end{corollary}
\proof
    Since $u$ is harmonic in every ball $B_r$ we can consider $r\to \infty$ in 
\autoref{theo:main-analysis} and use $\lim\limits_{r\to\infty}L(r)=0$ in 
order to prove that 
$u$ is a constant function.   
\qed

Our method to prove \autoref{lem:expans1} and \autoref{theo:main-analysis} is 
based on a purely analytic technique introduced in 
\cite{Sil06}. As mentioned above, a second aim of this work is to explain a 
probabilistic approach to results like \autoref{theo:main-analysis}. The 
starting point for these 
observations is that, for several linear differential or 
integrodifferential operators $A$, variants of \autoref{lem:growth} can be 
established with 
the help of corresponding Markov processes. Let
$X$ be the strong Markov process associated with the operator $A$, i.e. we
assume that the martingale problem has a unique solution. Denote by $T_A,
\tau_A$ the hitting resp. exit time for a measurable set $A \subset \R^d$
and by $\P_x$ the measure on the path space with $\P_x(X_0=x)=1$. The following
property then implies \autoref{lem:growth}.

\begin{proposition}\label{prop:hitting-prop}
There is a constant $c \in (0,1)$ such that for every $R>0$ and every
measurable set $A \subset B_{2R} \!\setminus\!B_R$ with $|(B_{2R} \!\setminus\!
B_R) \cap A| \geq \tfrac12 |B_{2R} \!\setminus\! B_R|$ and $x \in B_R$
\begin{align}\label{eq:hitting-estim}
\P_x (T_A < \tau_{B_{2R}}) \geq c \,.
\end{align}
\end{proposition}

This result is established for non-degenerate diffusions in \cite{KrSa79}, thus 
leading 
to a result like \autoref{theo:main-analysis} for elliptic differential 
operators of 
second order. The case of integrodifferential operators with fractional order 
of 
differentiability $\alpha \in
(0,2)$ is treated in \cite{BaLe02}. Therein it is shown that
\autoref{prop:hitting-prop} holds true for jump
processes $X$ generated by integral
operators $A$ of the form \eqref{eq:def_A} under the assumptions 
$K(x,h)=K(x,-h)$ and 
$K(x,h) \asymp |h|^{-d-\alpha}$ for
all $x$ and $h$ where
$\alpha \in (0,2)$ is fixed. Note that this class includes the case $\opA u =
-(-\Delta)^{\alpha/2} u$ and versions with bounded measurable coefficients.

\autoref{prop:hitting-prop} fails to hold for several cases that we are 
interested in. One 
example
is given by $\opA$ as in \eqref{eq:def_A} with 
$K(x,h)=k(h) \asymp |h|^{-d}$ for $|h|\leq 1$ and some appropriate condition for
$|h| > 1$. For example, the geometric stable process with its generator
$-\ln(1+(-\Delta)^{\alpha/2})$, $0 < \alpha \leq 2$, can be represented
by \eqref{eq:def_A} with a kernel $K(x,h)=k(h)$ with such a behavior for $|h|$
close to zero. The operator resp. the corresponding stochastic process can be
shown not to satisfy a uniformly
hitting estimate like \eqref{eq:hitting-estim}, see \cite{Mim14}. 

This leads 
to the question whether a priori estimates can be obtained by the approach from 
\cite{BaLe02} at all. In the second part of this work we address this question. 
It turns our that our main idea, i.e., to 
determine an new intrinsic
scale can also be used to establish a modification of 
\eqref{eq:hitting-estim}.
As we did in the proof of \autoref{lem:expans1}, we choose a measure 
different from the Lebesgue measure for the assumption
$|(B_{2R} \!\setminus\! B_R) \cap A| \geq \tfrac12 |B_{2R} \!\setminus\! B_R|$. 
We refer the reader to \autoref{sec:prob_estim} for further details

Since we employ methods from two different fields - methods from partial 
differential operators as in \cite{Sil06} and from stochastic analysis as in 
\cite{BaLe02} - it is interesting to compare both approaches. In both 
approaches we need to make several assumptions, e.g., solvability of the 
equation and existence of the corresponding Markov jump process. The 
conditions in the analysis approach are slightly less restrictive than 
those imposed in the approach using stochastic analysis. Note that, although we 
assume the solutions $u$ to be 
twice differentiable in the first part, the assertions resp. the 
constants in our results do not depend on the 
regularity of the functions $u$. Thus, the techniques and assertions presented 
here can be applied to nonlinear problems. 

\subsection{Examples} Let us look at different choices for the function $\ell$ 
used in condition \eqref{eq:A1}. Note that \eqref{eq:ell2} does not 
allow $\ell(h)$ to be zero (different from 
$K(x,h)$). Since the behavior of $\ell$ at zero is most important and 
characteristic, we provide examples of functions 
$\ell:(0,1) \to (0, \infty)$. For $a>1, r \in (0,1)$,  set 
$L(r)=\int_r^1 \ell (s)\frac{ds}{s}$ 
and $\varphi (r)=L ^{-1}(a^{-1} L (r))$.

\smallskip

\renewcommand{\arraystretch}{1.9}

\begin{table}[ht]
\caption{Different choices for a function $\ell$ when $\beta \in 
(0,2)$, $a > 1$.}
\begin{tabular}{ c  c  c  c  }
\hline
{\bf No.\! (i)} & $\ell_i(s)$ &
$L_i(s)$ & 
$\varphi_a(s)$ 
\\
\toprule
$1$ & $s^{-\beta} \, \ln(\frac{2}{s})^{2}$ & $\asymp s^{-\beta} \, \ln
(\frac{2}{s})^{2}$ & $\asymp 
s$\\
\hline
$2$ & $s^{-\beta}$ & $\frac{1}{\beta} (s^{-\beta}-1)$ & $\asymp 
s$\\
\hline
$3$ & $\ln(\frac{2}{s})$ & $\asymp \ln^2(\frac{2}{s})$ &
$\asymp s^{1/\sqrt{a}}$\\
\hline
$4$ & $1$ &  $\ln(\frac{1}{s})$ & $ s^{1/a}$\\
\hline
$5$ & $\ln(\frac{2}{s})^{-1}$ & $\asymp \ln(\ln(\frac{2}{s}))$ & $\asymp \exp(-
(\ln(\frac{2}{s}))^{1/a})$\\
\bottomrule
\end{tabular}
\label{tab:choices_l}
\end{table}

\subsection{Related results in the literature}
Let us comment on related results in the literature. The probabilistic 
approach, which we 
explain in \autoref{sec:intro-prob}, is based on the approach of \cite{BaLe02}.  
The analytic method, which we employ in the first part of the 
article, is based on \cite{Sil06}. Both approaches have been refined in many 
articles, allowing for more general kernels and treating fully nonlinear 
integrodifferential equations, but all these articles assume 
standard scaling properties, i.e., something like $K(x,h) 
\asymp |h|^{-d-\alpha}$ for some $\alpha \in (0,2)$. We refer to the 
discussions in \cite{KaSc14, SiSc14} for further references. Note that our 
regularity result \autoref{theo:main-analysis} is stronger than 
\autoref{theo:main-prob} because we can allow for right-hand sides $f$ in the 
integrodifferential equation and for more general kernels.  

The current work comprises the two preprints \cite{KaMi13} and \cite{KaMi14} 
where the approaches by analytic and probabilistic methods are explained 
separately. After \cite{KaMi13} had appeared, several articles have made use of 
the ideas therein. In 
\cite{Bae14} nonlocal problems are studied where the kernels are supposed to 
satisfy certain upper and lower scaling conditions. These 
assumptions do not include limit cases like \eqref{eq:exa_L} since 
some comparability with kernels like $|h|^{-d-\alpha}$ for 
$\alpha \in (0,2)$ is still assumed. In \cite{KKL14} the authors study fully 
nonlinear problems with similar assumptions on the kernels as in \cite{Bae14}. 
In \cite{ChZh14} the authors extend the regularity estimates of \cite{KaMi13} to 
time-dependent equations with drifts. The article \cite{JaWe14} is not directly 
related to \cite{KaMi13} but mentions the need to consider $f \ne 
0$. We 
solve this problem.

\subsection{Organization of the article} In \autoref{sec:levy} we review
the relation between translation invariant nonlocal operators and
semigroups/L\'{e}vy processes. Presumably, \autoref{prop:relation} is of some 
interest to many readers since it establishes a
one-to-one relation between the behavior of a L\'{e}vy measure at zero and the
multiplier of the corresponding generator for large values of $|\xi|$. 
\autoref{sec:proof-expans1} and \autoref{sec:proof-theo-analysis} contain the 
proof of \autoref{lem:expans1} and \autoref{theo:main-analysis} respectively. 
In \autoref{sec:intro-prob} we explain the probabilistic approach to 
\autoref{theo:main-analysis}, which leads to \autoref{theo:main-prob}. Note that 
we are slightly changing the assumptions there. The probabilistic approach is 
based on a Krylov-Safonov type hitting lemma, which is 
\autoref{prop:hitting_new}. \autoref{sec:prob_estim} 
contains the proof of this result and of 
\autoref{theo:main-prob}. The last section is 
\autoref{sec:appendix} in
which we collect important properties of regularly resp. slowly varying
functions. 

\subsection{Acknowledgement} The authors thank Tomasz Grzywny and Jongchun 
Bae for discussions about the assumptions used in this article. Further thanks 
are due to an anonymous referee for comments which helped us to improve 
the presentation.

\section{Multipliers and L\'{e}vy measure: Analysis meets 
Probability}\label{sec:levy}

The aim of this section is to provide some background about translation 
invariant integrodifferential operators and related stochastic processes. The 
results explained here motivated the search for a new scale function which is 
a key element of the whole project. However, the material of this section is 
not needed for the proofs of the main results. 

In this paper we provide two approaches to \autoref{theo:main-analysis}. One 
approach uses techniques from analysis, the other one uses stochastic 
processes. Note that the quantity $K(x,h) \d h$ in \eqref{eq:def_A} has a clear 
interpretation in terms of probability. For fixed $x$, the quantity $\int_M 
K(x,h) \d h$ describes the intensity with which the corresponding process 
performs jumps from some point $x \in \R^d$ to a point from the set $x + M$. 
In this 
sense, the conditions 
\eqref{eq:ell1}-\eqref{eq:ell3} say something about the behavior of the process. 
On the other hand, the conditions say something about mapping properties of the 
operator $\opA$. In this section we explain the link between these two 
viewpoints. We restrict ourselves to the cases of translation invariant 
operators, i.e., we assume $K(x,h)$ to be independent of the first variable. 
This allows us to give a focused presentation. Note that the results of the 
section are not used in the rest of the article. 

In the translation invariant case, i.e. when 
$K(x,h)$ does not depend on $x$, there is a one-to-one
correspondence between $\opA$ and multipliers, semigroups and stochastic
processes. One aim is to prove how the behavior of $\ell(|h|)$ for small values
of $|h|$ translates into properties of the multiplier or characteristic
exponent $\psi(|\xi|)$ for large values of $|\xi|$. This is achieved in
\autoref{prop:relation}. We add a subsection where we discuss which
regularity results are known in critical cases of the (much simpler) translation
invariant case. Note that our set-up, although
allowing for an irregular dependence of $K(x,h)$ on $x \in \R^d$, leads to new
results in these critical cases. 

\subsection{Semigroups, generators and L\'{e}vy processes}

A stochastic process $X=(X_t)_{t\geq 0}$ on a probability space
$(\Omega,\mathcal{F},\P)$ is called a L\' evy process if it has
stationary and independent increments, $\P(X_0=0)=1$  and its
paths are $\P$-a.s. right continuous with left limits. For $x\in \R^d$ we
define a $\P_x$
to be the law of the process $X+x$. In particular, $\P_x(X_t\in B)=\P(X_t\in
B-x)$ for $t\geq 0$ and measurable sets $B\subset \R^d$.

Due to stationarity and independence of increments, the characteristic exponent 
of $X_t$ is given by
\[
 \E[e^{i\langle\xi, X_t\rangle}]=e^{-t\psi(\xi)},
\]
where $\psi$ is called characteristic exponent of $X$. It has the
following L\'
evy-Khintchine representation
\begin{equation}\label{eq:lk}
 \psi(\xi)=\frac{1}{2}\langle A\xi,\xi\rangle+\langle b,
\xi\rangle+\int_{\R^d\setminus
\{0\}}(1-e^{i\langle \xi,
h\rangle}+i\langle \xi,h\rangle \mathbbm{1}_{B_1}(h))\nu(dh)\,,
\end{equation}
where $A$ is a symmetric non-negative definite matrix, $b\in \R^d$ and $\nu$ is
a measure on
$\R^d\setminus \{0\}$ satisfying $\int_{\R^d\setminus\{0\}}(1\wedge
|y|^2)\nu(dy)<\infty$ called the
L\' evy measure of $X$. 

The converse also holds; that is, given $\psi$ as in the  L\' evy-Khintchine
representation \eqref{eq:lk}, there exists a L\' evy process $X=(X_t)_{t\geq
0}$ with the characteristic exponent $\psi$. The equality \eqref{eq:lk} 
provides a link to an analytic viewpoint on L\'evy processes. If $\nu$ is a 
L\'evy measure, i.e., a Borel measure on $\cB(\R^d \setminus \{0\})$, then one 
can construct a convolution semigroup of probability measures $(\nu_t)_{t > 0}$ 
such that the Fourier transform of $\nu_t$ equals $e^{-t \psi}$ with $\psi$ as 
in \eqref{eq:lk}. This approach can be found in \cite{BeFo75}. 

Let $X=(X_t)_{t\geq 0}$ be a L\' evy process
corresponding to the characteristic exponent $\psi$ as in \eqref{eq:lk} with
$A=0$, $b=0$ and a L\' evy measure $\nu(dh)$. Then $P_t 
f(x):=\E_x[f(X_t)]$ defines a strongly continuous contraction semigroup
of operators $(P_t)_{t\geq 0}$ on the space of bounded 
uniformly continuous functions on $\R^d$ equipped with
the supremum norm.  Moreover, it is a
convolution semigroup, since 
\[
 \P_t f(x)=\E_0[f(x+X_t)]=\int_{\R^d}f(x+y)\nu_t(dy)\,,
\]
with $\nu_t(B):=\P(X_t\in B)$. The infinitesimal generator ${\opA}$ of the 
semigroup $(P_t)_{t\geq 0}$ is
given by
\begin{align}\label{eq:generator_levy}
 {\opA}u(x)=\int_{\R^d\setminus\{0\}}\big(u(x+h)-u(x)-\langle \nabla
u(x),h\rangle \mathbbm{1}_{B_1}(h)\big) \nu(\!\d h) \,,
\end{align}
if $u$ is sufficiently regular, see the proof of \cite[Theorem 
31.5]{Sat99}. Note that the process $\big(u(X_t)-u(X_0)-\int_0^t{\opA}u(X_s)\d 
s)_{t \geq 0}$ is a
martingale (with respect to the natural filtration) for every $u\in 
C_b^2(\R^d)$, see the proof of \cite[Proposition VII.1.6]{ReYo99}. In this 
sense the process $X$ corresponds to the given L\'evy measure $\nu$ and, in our 
set-up, to the kernel $K(x,h)=k(h)$. For details about L\' evy processes we 
refer to \cite{Ber96,Sat99}.

Let us now explain the connection between 
the  characteristic exponent $\psi$ and the symbol of the operator ${\opA}$. To 
be more precise, if
$\hat{f}(\xi)=\int_{\R^d}e^{i\langle \xi, x\rangle} f(x)\d x$ denotes the 
Fourier transform
of a function $f\in L^1(\R^d)$, then 
\[
 \widehat{{\opA}f}(\xi)=-\psi(-\xi)\hat{f}(\xi)
\]
for any $f\in \mathcal{S}(\R^d)$, where $\mathcal{S}(\R^d)$ is the Schwartz
space, see \cite[Proposition I.2.9]{Ber96}. Hence, $-\psi(-\xi)$ is the
symbol (multiplier) of the operator ${\opA}$. The following result explains 
how, in the case $\nu(\d h)=k(h) \d h = K(x,h) \d h$, the  kernel $K(x,h)=k(h)$ 
is related to the characteristic exponent resp. the multiplier.

\begin{proposition}\label{prop:relation} Assume that the operator $\opA$ 
defined 
on $\cS$ 
is given by
\eqref{eq:generator_levy}. Assume $\nu(\!\d h)=k(h) \d h = K(x,h) \d 
h$ where $K\colon \R^d\times(\R^d\setminus\{0\})\rightarrow 
[0,+\infty)$ is a measurable function with $K(x,h)=K(x,-h)$ for almost all 
$x,h$. We assume that $K$ satisfies \eqref{eq:A1}, \eqref{eq:K-Levy-bound} with 
$R_0 \in (0,+\infty]$. Set 
$L(r)=\int_r^{R_0} \ell(s)\frac{ds}{s}$. Then there 
are constants $c>0$ and $r_0>0$ such that 
\[
 c^{-1}L(|\xi|^{-1})\leq \psi(\xi)\leq cL(|\xi|^{-1}) \quad \text{ for } \xi\in
\R^d,\
|\xi|\geq r_0\,.
\]
\end{proposition}

The assumptions of \autoref{prop:relation} allow to treat sophisticated 
examples. However, it is instructive to think about the simple examples 
\begin{align*}
K(x,h)&= |h|^{-d-\alpha} \text{ for some } \alpha \in (0,2) \,, \\
K(x,h)&= |h|^{-d} \mathbbm{1}_{B_1}(h) \,, \\
K(x,h)&= |h|^{-d} \ln(2/|h|)^{\pm 1} \mathbbm{1}_{B_1}(h) \,.
\end{align*}

\proof
Note first that, by \eqref{eq:A1},
\begin{align*}
 \kappa^{-1}j(|h|)\leq k(h)\leq \kappa j(|h|) \text{ for } 0 < |h| < 
R_0 
\end{align*}
for $j(s) = s^{-d} \ell(s)$. Since $1-\cos{x}\leq \frac{1}{2}x^2$, it follows 
from \eqref{eq:ell2} and \autoref{lem:L} that  
\begin{align*}
 \psi(\xi)&\leq \tfrac{1}{2}|\xi|^{2}\int_{|h|\leq
|\xi|^{-1}}|h|^2j(|h|)\d h+2\int_{ |\xi|^{-1}<|h|<R_0}j(|h|)\d h+2\int_{|h|\geq R_0}k(h)\,dh\\
&\leq 
c_1\left[|\xi|^2\int_0^{|\xi|^{-1}}s\ell(s)\d s+L(|\xi|^{-1})+1\right]\\&\leq 
c_2\left[|\xi|^2\ell(|\xi|^{-1})\int_0^{|\xi|^{-1}}sc_L^{-1}(\tfrac{s}{|\xi|^{-1
}})^{
-\gamma}\d s+L(|\xi|^{-1})\right]\\
&\leq c_3 (\ell(|\xi|^{-1})+L(|\xi|^{-1}))\leq c_4 L(|\xi|^{-1})\,.
\end{align*}

In order to prove the lower bound, we employ an idea of \cite{Grz14}. Let us 
first consider the case $R_0 = + \infty$. We choose an orthogonal 
transformation of the
form $Oe_1=|\xi|^{-1}\xi$,
where $e_1:=(1,0,\ldots,0)\in \R^d$. Then a change of variable 
yields 
\begin{align*}
 \psi(\xi)&=\int_{\R^d\setminus\{0\}}(1-\cos( \langle \xi,
h \rangle))j(|h|)\d h=\int_{\R^d\setminus\{0\}}(1-\cos{(|\xi|h_1)})j(|h|)\d h\,.
\end{align*}
By the Fubini theorem,
we obtain
\[
\psi(\xi)\geq 2\int_0^\infty (1-\cos{(|\xi|r)})F(r)\d r,
\]
where $
F(r):=\int_{\R^{d-1}}j(\sqrt{|z|^2+r^2})\d z,\quad r>0$.
Using \eqref{eq:ell3} we deduce for every $0<r\leq s$, 
\[
    F(r)=\int_{\R^{d-1}} 
\left(\frac{|z|^2+s^2}{|z|^2+r^2}\right)^{d/2}\frac{\ell(\sqrt{|z|^2+r^2})}{
\ell(\sqrt{
|z|^2+s^2})}j(\sqrt{|z|^2+s^2})\d z\geq c_U^{-1}F(s)\,.
\]

Now, 
\begin{align*}
 \psi(\xi)&\geq 2 \sum_{k=0}^\infty 
\int_{|\xi|^{-1}(\frac{\pi}{2}+2k\pi)}^{|\xi|^{-1}(\frac{3\pi}{2}+2k\pi)}(1-\cos
{(|\xi|r)})F(r)\,
dr
\geq \frac{c_U^{-1}\pi}{|\xi|}\sum_{k=0}^\infty 
F(|\xi|^{-1}(\tfrac{3\pi}{2}+2k\pi))\\
&\geq c_U^{-2}\sum_{k=0}^\infty 
\int_{|\xi|^{-1}(\frac{3\pi}{2}+2k\pi)}^{|\xi|^{-1}(\frac{3\pi}{2}+(2k+1)\pi)}
F(r)\d r
\geq c_U^{-2}\int_{\frac{3\pi}{2}|\xi|^{-1}}^\infty  F(r)\d r\\
&\geq c_{5} \int_{|h|\geq \frac{3\pi}{2}|\xi|^{-1}} j(|h|)\d h=c_{6}
L(\tfrac{3\pi}{2}|\xi|^{-1})\geq c_{7} L(|\xi|^{-1})\,,
\end{align*}
where in the last inequality we have used \autoref{lem:L}. 
The case $R_0<+\infty$ can be proved similarly.
\qed

\subsection{Some related results from potential theory}
Let us explain which results, related to \autoref{theo:main-analysis}, have been
obtained in the case where $K(x,h)$ is independent of $x \in \R^d$. In this 
case, methods from potential theory can be used. 

H\"{o}lder estimates of harmonic functions are obtained for the L\' evy process
with the characteristic exponent
$\psi(\xi)=\frac{|\xi|^2}{\ln(1+|\xi|^2)}-1$ in \cite{Mim13} by establishing
a Krylov-Safonov type estimate replacing the Lebesgue measure with the capacity
of the sets involved. Recently, regularity estimates have been
obtained in \cite{Grz14} for a class of isotropic unimodal L\' evy processes 
which
is quite general but does not include L\' evy processes with slowly varying L\'
evy exponents such as geometric stable processes. Regularity of harmonic
functions for such processes is investigated in \cite{Mim14}, where it is shown
that a result like \autoref{prop:hitting-prop} fails. Using the 
Green
function, logarithmic bounds for the modulus of continuity are obtained. At this
point it is worth mentioning that the transition density $p_t(x,y)$ of the 
geometric stable process satisfies $p_1(x,x)=\infty$, see \cite{SSV06}. This
illustrates that regularity results like \autoref{theo:main-analysis} in the 
case $K(x,h)= |h|^{-d} \mathbbm{1}_{B_1}(h)$ and in similar cases are quite 
delicate.

\section{Proof of \autoref{lem:expans1}}\label{sec:proof-expans1}

Before we proceed to the main proof, let us provide two auxiliary statements 
which 
indicate the link of the scale function $\varphi$ with the kernels $K$.
\begin{lemma}\label{lem:L} The following properties of the function $L$ hold 
true: 
\begin{alignat}{2}
&L(r)\geq \gamma^{-1}c_L \,\ell(r)\quad &&\text{ for } r>0 \, \tag{i}\\
\text{ In the case } R_0 = + \infty: \quad  &\frac{L(r \lambda)}{L(r)} \geq 
c_L \lambda^{-\gamma} \qquad &&\text{ for } \lambda \geq 1, r>0\,, \tag{ii-a} \\
   \text{ In the case } R_0 < + \infty: \quad &\frac{L(r \lambda)}{L(r)}\geq 
\frac{c_L}{2} \lambda^{-\gamma} &&\text{ for } 1 \leq \lambda < \lambda_1, 0 
< r < r_1\,, \tag{ii-b}
\end{alignat}
where $\lambda_1 > 1$ can be chosen arbitrarily and $r_1 > 0$ 
depends on $\lambda_1$ and $R_0$.
\end{lemma}

\proof
  From (\ref{eq:ell2}) we deduce the following
  \[
   L(r)=\ell(r)\int_r^{R_0}\tfrac{\ell(s)}{\ell(r)}\tfrac{ds}{s}\geq c_L 
\ell(r) r^{\gamma}\int_r^{R_0} s^{-1-\gamma}\d s=c_L\gamma^{-1}\ell(r)\,,
  \]
which proves part (i). In the case $R_0 = + \infty$ we obtain
  \[
   L(r \lambda)=\int_{r \lambda}^\infty \ell(s)\tfrac{ds}{s}=\int_r^\infty 
\ell(s \lambda)\tfrac{ds}{s}=\int_r^\infty \tfrac{\ell(s 
\lambda)}{\ell(s)}\ell(s)\tfrac{ds}{s}\geq c_L \lambda^{-\gamma}\int_r^\infty 
\ell(s) \tfrac{ds}{s}=c_L \lambda^{-\gamma}L(r)\,,
  \]
which is one part of claim (ii). In the case $R_0 < +\infty$ we deduce 
  \begin{align*}
   L(r \lambda) &=\int_{r \lambda}^{R_0} \ell(s)\tfrac{ds}{s} = 
\int_r^{R_0/\lambda} \ell(s \lambda)\tfrac{ds}{s} \geq c_L \lambda^{-\gamma} 
\int_r^{R_0/\lambda} \ell(s)\tfrac{ds}{s}  \geq c_L \lambda^{-\gamma} 
\big( L(r) - L(R_0/\lambda) \big) \\
& \geq  c_L \lambda^{-\gamma} 
\big( L(r) - L(R_0/\lambda_1) \big)  =  \frac{c_L}{2} \lambda^{-\gamma} 
L(r) + \frac{c_L}{2} \lambda^{-\gamma} \big( L(r) -  2 L(R_0/\lambda_1) 
\big) \geq \frac{c_L}{2} \lambda^{-\gamma} L(r)  \,, 
  \end{align*}
which proves the remaining case and completes the proof.   
\qed

\begin{lemma}\label{lem:phi}
 Assume $a>1$. Then $\mu(B_{\varphi_a(r)}\setminus B_r)=|\partial B_1|\ln a$,
 where $|\partial B_1|$ denotes the surface area of the unit sphere in $\R^d$.
\end{lemma}
\proof
The proofs follows by introducing polar coordinates: 
  \[
   \mu(B_{\varphi_a(r)}\setminus B_r)=|\partial 
B_1|\int_r^{\varphi_a(r)}\tfrac{1}{L(s)}\tfrac{\ell(s)\d s}{s}=|\partial 
B_1|\ln\tfrac{L(r)}{L(\varphi_a(r))}=|\partial B_1|\ln a\,.
  \]
\qed

\begin{lemma}\label{lem:j-func}
Set $j(s) = s^{-d} \ell(s)$ for $0<s<R_0$. Let $M \geq 1$. Then
\begin{align*}
 s \leq Mt \leq R_0
 \quad \text{ and } \quad t\leq R_0  \quad \text{ imply } \quad j(t) \leq c j(s) 
\,,
\end{align*}

with $c= \max\{c_U, M^{\gamma+d} c_L^{-1}\}$.
\end{lemma}
\proof
Assume $s,t >0$ with $s \leq Mt \leq R_0$ for some $M \geq 1$. We consider two 
cases. If $s < t$, then 
\[ j(t) = t^{-d} \ell(t) = (s/t)^d s^{-d} \ell(s (t/s)) \leq c_U s^{-d} 
\ell(s) \,, \]
where we have applied \eqref{eq:ell3}. If $t \leq s \leq Mt$, then
\[ j(t) = t^{-d} \ell(t) \leq (s/t)^{\gamma} c_L^{-1} (M^{-1}s)^{-d} \ell(s) \leq 
M^{\gamma+d} c_L^{-1} j(s)  \,, \]
where we have applied \eqref{eq:ell2}. The proof is complete.
\qed

We are now able to provide the proof of our first main result. Recall that this 
result is proved under the assumption $R_0 = + \infty$.

\proof[Proof of \autoref{lem:expans1}] 
Define $\beta\colon [0,+\infty)\rightarrow [0,+\infty)$, $\beta (r) = 
\exp(-r^2)$, and further 
  \[
 b(x):=\beta(|x|)\quad \text{ and }\quad 
b_r(x):=\beta_r(x):=\beta(r^{-1}|x|)\quad \text{ for }\  x\in \R^d, r>0\,.
  \]
First we estimate $-A 
b_r$. For $r > 0$ we can deduce from  
(\ref{eq:ell1}), (\ref{eq:ell2}) and \autoref{lem:L} the following:
\begin{align*}
   -A b_r(x)&= \int_{\R^d}\left(b(\tfrac{x}{r})-b(\tfrac{x+y}{r}) +
   \langle \nabla b(\tfrac{x}{r}) , \tfrac{y}{r} \rangle 
\mathbbm{1}_{B_1}(y)\right) K(x,y)
\d y \\
   & \leq 
c_1\int_{\R^d}\left(\left(\tfrac{|y|}{r}\right)^2 
\mathbbm{1}_{B_r}(y)+ \mathbbm{1}_{B_r^c}
(y)\right)\tfrac{
\ell(|y|)}{|y|^d}\d y\\
   & =c_2\left(r^{-2}\int_0^{r}s\ell(s)\d s+\int_{r}^{\infty} 
\ell(s)\tfrac{ds}{s} \right)\\
   & =c_2 \left(r^{-2}\ell(r)\int_0^{r}s\tfrac{\ell(s)}{\ell(r)}\d 
s+L(r) \right)\\
   &\leq c_3 \left(r^{-2}\ell(r)\int_0^{r}s\left(\tfrac{r}{s}\right)^{\gamma}\d 
s+
 L(r) \right)\\
   &\leq c_4 (\ell(r)+L(r) )\leq c_5 L(r)\,,
  \end{align*}
Hence, 
  \begin{equation}\label{eq:b_r}
   \sup_{x\in \R^d}-A b_r(x)\leq c_6 L(r)\,.
  \end{equation}

  Set $ \theta:=\frac{1}{a}(\beta(1)-\beta(\frac{3}{2}))=\frac{1}{a}(\beta_r( 
r)-\beta_r(\frac{3r}{ 2}))$, where $a >2$ will be chosen later 
independently of $v$ and $r$.
  
We claim that one can choose $a>2$ so large that $v(x)\leq 1- \theta$ for $x\in B_r$. 
Assume that this is not true. Then for $a>2$ there is $x_0\in B_r$ satisfying
  \[
   v(x_0)\geq 1-
    \theta=1-a^{-1}\beta_r(r)+a^{-1}\beta_r(\tfrac{3r}{2})\,.
  \]
  Since $|x_0|<r$, 
  \begin{align}\label{eq:tmp123}
  \begin{split}
   v(x_0)+a^{-1}b_r(x_0)&\geq 1+a^{-1}\beta_r(|x_0|)-a^{-1}\beta_r(r)+a^{-1}\beta_r(\tfrac{3r}{2}) \\
   &>1+a^{-1}\beta_r(\tfrac{3r}{2})\\
   &\geq v(y)+a^{-1}\beta_r(|y|)\quad \text{ for all }\quad y\in B_{\varphi(r)}\setminus B_{\frac{3r}{2}},
  \end{split}
  \end{align}
  where the last inequality follows from the assumption $v(x)\leq 1$ for $x\in 
B_{\varphi(r)}$ and $\beta_r(|y|)\leq \beta_r(\frac{3r}{2})$. By choosing $a$ 
sufficiently large, we will make sure that $\varphi(r) > \frac{3r}{2}$. It 
follows from \eqref{eq:tmp123} that $v+a^{-1}b_r$ attains its maximum at $x_1\in 
B_{\frac{3r}{2}}$ and $(v+a^{-1}b_r)(x_1)\geq 1+a^{-1}\beta_r(\tfrac{3r}{2}) 
>1$.

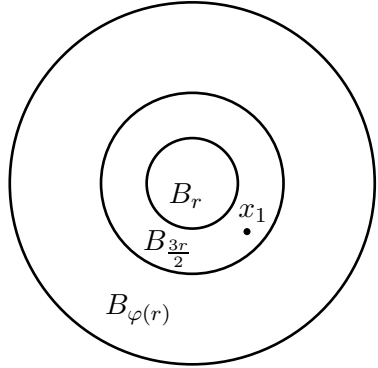
\begin{figure}[h]
\begin{tikzpicture}[scale=0.8]
\draw[line width =1pt] (0, 0) circle (0.75);
\draw[line width =1pt] (0, 0) circle (1.5);
\draw[line width =1pt] (0, 0) circle (3);
\draw[fill] (.9,-0.8) circle [radius=1.4pt];
\node[above] at (1,-0.8) {$x_1$};
\node[right] at (-1.6,-2.1) {$B_{\varphi(r)}$};
(-1.06,-1.06);
\node[below] at (-0.4,-0.6) {$B_{\frac{3r}{2}}$};
\node at (-0.1,-0.2) {$B_r$};
\end{tikzpicture}
\caption{$B_r \subset B_{\frac{3r}{2}} \subset B_{\varphi(r)}$} 
\label{fig:circles}
\end{figure}

  The idea now is to establish a contradiction by evaluating $-A(v+a^{-1}b_r)(x_1)$ in two different ways. First, by (\ref{eq:b_r})
  \[
   -A(v+a^{-1}b_r)(x_1)\leq a^{-1}L(r)+c_6a^{-1}L(r)=(1+c_6)a^{-1}L(r)\,.
  \]

   On the other hand, since $v+a^{-1}b_r$ attains its maximum at $x_1$, 
   $\nabla(v+a^{-1}b_r)(x_1)=0$ and hence
   \begin{align*}
    -A(v+a^{-1}b_r)(x_1)&
    = \int_{\R^d\setminus 
\{0\}}\left((v+a^{-1}b_r)(x_1)-(v+a^{-1}b_r)(x_1+y)\right)K(x_1,y)\d y\\
	      &=\int_{\{y\in \R^d\setminus\{0\}\colon x_1+y\in 
B_{\varphi(r)}\}}\left((v+a^{-1}b_r)(x_1)-(v+a^{-1}b_r)(x_1+y)\right)K(x_1,y)\d 
y\\
	      &\ +\int_{\{y\in \R^d\setminus\{0\}\colon x_1+y\not\in 
B_{\varphi(r)}\}}\left((v+a^{-1}b_r)(x_1)-(v+a^{-1}b_r)(x_1+y)\right)K(x_1,y)\d 
y\\
	      &=:I_1+I_2\,.
   \end{align*}

Since $(v+a^{-1}b_r)$ attains its maximum on $B_{\varphi(r)}$ in $x_1$ with 
$(v+a^{-1}b_r) >1$,
    by (\ref{eq:A1}) we obtain
    \begin{align*}
     I_1&\geq \int_{\{y\in \R^d\setminus\{0\}\colon x_1+y\in 
B_{\varphi(r)}\setminus B_r,\ v(x_1+y)\leq 
0\}}\left((v+a^{-1}b_r)(x_1)-(v+a^{-1}b_r)(x_1+y)\right)K(x_1,y)\d y\\
     &\geq c_7 (1-a^{-1}\|b\|_\infty) \int_{\{y\in \R^d\setminus\{0\}\colon 
x_1+y\in B_{\varphi(r)}\setminus B_r,\ v(x_1+y)\leq 0\}} j(|y|)\d y\,,
    \end{align*}
  with $j(s):=s^{-d}\ell(s)$.
  
  Using $|y|\leq |x_1+y|+|x_1|\leq |x_1+y|+\frac{3r}{2}\leq \frac{5}{2}|x_1+y|$ for 
$x_1+y\in B_{\varphi(r)}\setminus B_r$,  we deduce from (\ref{eq:ell2}) that $j(|y|)\geq 
c_8 j(|x_1+y|)$. Here we have applied \autoref{lem:j-func}. The 
assumptions of the lemma imply 
  \begin{align*}
   I_1&\geq c_{9} (1-a^{-1}\|b\|_\infty) \int_{\{y\in \R^d\setminus\{0\}\colon 
x_1+y\in B_{\varphi(r)}\setminus B_r,\ v(x_1+y)\leq 0\}} j(|x_1+y|)\d y\\
      &=c_{9} (1-a^{-1}) \int_{\{y\in \R^d\setminus\{0\}\colon y\in 
B_{\varphi(r)}\setminus B_r,\ v(y)\leq 0\}} j(|y|)\d y\\
      & =c_{9}(1-a^{-1})\int_{\{y\in \R^d\setminus\{0\}\colon y\in B_{\varphi(r)}\setminus 
B_r,\ v(y)\leq 0\}} L(|y|)\,\mu(dy)\\
      &\geq c_{9}(1-a^{-1})L(\varphi(r))\mu((B_{\varphi(r)}\setminus B_r)\cap \{v\leq 
0\})\\
      &\geq c_{9}(1-a^{-1})a^{-1}L(r)\delta \mu(B_{\varphi(r)}\setminus B_r)\\
      &\geq c_{10}(1-a^{-1})\delta L(r)\tfrac{\ln a}{a},
  \end{align*}
where in the last inequality we have used \autoref{lem:phi}.

\autoref{lem:L} implies that, if we consider $a>c_L^{-1}(5/2)^{\gamma}$, then  
$\frac{L(r)}{L(\frac{5r}{2})}\leq c_L^{-1}(5/2)^{\gamma}$. Hence
\[
 L(\varphi(r))=a^{-1}\frac{L(r)}{L(\frac{5r}{2})}L(\tfrac{5r}{2})\leq L(\tfrac{5r}{2})
\]
and, since $L$ is decreasing, we obtain that $\varphi(r)\geq \frac{5r}{2}$.
To estimate $I_2$ we note that for $x_1+y\not\in B_{\varphi(r)}$ it follows 
from (\ref{eq:tmp123}) 
that 
\[
(v+a^{-1}b_r)(x_1)-a^{-1}b_r(x_1+y)\geq 1+a^{-1}\beta_r(\tfrac{3r}{2})-a^{-1}\beta_r(\varphi(r))\geq 1\,.
\]
Hence, this together with the growth assumption on $v$ yields
\begin{align*}
 I_2&\geq -c_{11}\int_{\{x_1+y \in \R^d \setminus 
B_{\varphi(r)}\}}\left(\tfrac{L(\varphi(r))}{L(|x_1+y|)}\right)^{\eta}j(|y|)\d 
y\,.
\end{align*}

Note that $x_1+y\not\in B_{\varphi(r)}$ implies $|y|\geq |x_1+y|-|x_1|\geq 
\tfrac{5r}{2}-\tfrac{3r}{2}=r$ and  
\[
	|x_1+y|\leq \tfrac{3r}{2}+|y|\leq \tfrac{3}{2}|y|+|y|=\tfrac{5}{2}|y| \text{ and 
} y\not\in B_{\frac{2}{5}\varphi(r)}\,.
\]
In this case $L(|x_1+y|)\geq L(\frac{5}{2}|y|)$ and \[{\{y\in \R^d\setminus\{0\}\colon 
x_1+y\not\in B_{\varphi(r)}\}}\subset {\{y\in \R^d\setminus\{0\}\colon |y|\geq 
\tfrac{2}{5}\varphi(r)\}} \,.\]

Thus we obtain
\begin{align*}
 I_2&\geq -c_{12} \int_{\{x_1+y \in \R^d \setminus 
 B_{\varphi(r)}\}}\left(\tfrac{L(\varphi(r))}{L(\frac{5}{2}|y|)}\right)^\eta 
j(|y|)\d y \\
 &\geq -c_{13}\int_{\frac{2}{5}\varphi(r)}^{\infty} 
\left(\tfrac{L(\varphi(r))}{L(\frac{5}{2}s)}\right)^\eta\,\tfrac{\ell(s)}{s}ds\\
 &\geq -c_{14}\int_{\frac{2}{5}\varphi(r)}^{\infty} 
\left(\tfrac{L(\varphi(r))}{L(s)}\right)^\eta\,(-L'(s))ds \\
 &=-c_{14}L(\varphi(r))^\eta(1-\eta)^{-1}L(\tfrac{2}{5}\varphi(r))^{1-\eta}\\
 &\geq -c_{15}L(\varphi(r)) =-c_{16}\frac{L(r)}{a}\,,
\end{align*}
where in the third inequality \autoref{lem:L} and in the last inequality 
monotonicity of $L$ and \autoref{lem:L} again have been used.

Finally, we obtain
\[
 (1+c_6)\tfrac{L(r)}{a}\geq c_{10}(1-a^{-1})\delta L(r)\tfrac{\ln 
a}{a}-c_{16}\tfrac{L(r)}{a}\,,
\]
or written in another way $(1+c_6+c_{16}) \geq 
c_{10}(1-a^{-1})\delta \ln a$.
Choosing $a>2$ large enough leads to a contradiction. 

This means that we have proved that there exists $a>2$ such that \[v(x)\leq 1-a^{-1}(\beta(1)-\beta(\tfrac{3}{2}))=1-\theta\quad \text{for all }\quad x\in B_r\,.\]

Note that our choice of $a$ does not depend on $r$; hence the assertion of the 
theorem holds for every $r>0$ with the same choice of $a$ and $\theta$.

\qed

\section{Proof of \autoref{theo:main-analysis}}\label{sec:proof-theo-analysis}

\proof
First of all, we restrict the general case and assume that 
conditions \eqref{eq:A1} and \eqref{eq:K-Levy-bound} hold true with $R_0 = + 
\infty$. At the end of the proof we explain how to reduce the general case to 
this case. 

Let $r > 0$. Assume $u \in 
C^2_b(\R^d)$  satisfies $Au=f$ in $B_r$ where $f$ is essentially bounded. We assume $u 
\not \equiv 0$ and prove assertion \eqref{eq:theo-estim} in the simplified 
case 
$\|u\|_\infty 
\leq 1/2$ and $\|f\|_{L^\infty(B_r)} \leq \tfrac12 L(r/2)$. Let us briefly 
explain why this 
is 
sufficient. In the general case we would set 
\[ \widetilde{u} = \frac{u}{2 \|u\|_\infty + 2 L(r/2)^{-1} \|f\|_{L^\infty(B_r)}} \]  
If $u$ solved $Au=f$ in $B_r$, then $\widetilde u$ would solve $A 
\widetilde{u}=\widetilde{f}$ in $B_r$ with $\|\widetilde{u}\|_\infty \leq 1/2$ and 
$\|\widetilde{f}\|_{L^\infty(B_r)} \leq \tfrac12 L(r/2)$. Thus we could apply the result 
in the simplified case and obtain:
 \begin{align}\label{eq:theo-estim-modified}
  \sup_{x,y\in B_{r/4}}\tfrac{|u(x)-u(y)|}{L(|x-y|)^{-\beta}} 
  & \leq \left( \tfrac{c}{2} L(r)^{\beta} + \tfrac{c}{2} L(r/2) L(r)^{\beta-1} \right) \; 
\left(2 \|u\|_\infty + 
2 L(r/2)^{-1} \|f\|_{L^\infty(B_r)} \right) \\
 & \leq \widetilde{c} L(r)^{\beta} \|u\|_\infty +  \widetilde{c}  
L(r)^{\beta-1} \|f\|_{L^\infty(B_r)}  
\,,
\end{align}
where $\widetilde{c}$ is another constant, depending on $c_L$ and $\gamma$ because 
of \autoref{lem:L}.

Hence we can restrict ourselves to $\|u\|_\infty \leq 
1/2$ and $\|f\|_{L^\infty(B_r)} \leq \tfrac12 L(r/2)$. Let $x_0\in B_{r/4}$. 
Without loss of generality we may assume $u(x_0) > 0$.

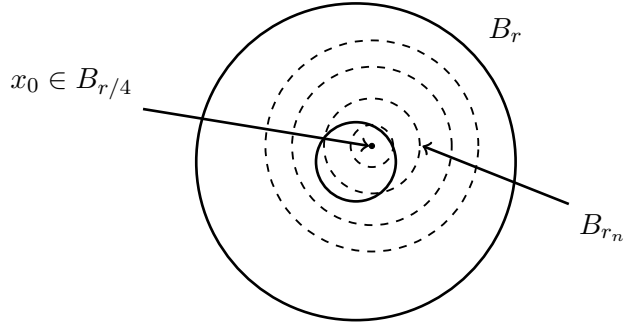
\begin{figure}[h]
\begin{tikzpicture}[scale=0.7]
\draw[line width =1pt] (0, 0) circle (0.75);
\draw[fill] (0.3,0.3) circle [radius=1.4pt];
\draw[line width =1pt] (0, 0) circle (3);
\draw[line width =0.7pt, dashed] (0.3, 0.3) circle (2);
\draw[line width =0.7pt, dashed] (0.3, 0.3) circle (1.5);
\draw[line width =0.7pt, dashed] (0.3, 0.3) circle (0.9);
\draw[line width =0.7pt, dashed] (0.3, 0.3) circle (0.4);
\node[right] at (2.3,2.5) {$B_r$};
\draw [->, line width =1pt] (-4,1) -- (0.25,0.3);
\node[above left] at (-4,1) {$x_0 \in B_{r/4}$};
\draw [->, line width =1pt] (4,-0.8) -- (1.25,0.3);
\node[below right] at (4,-0.8) {$B_{r_n}$};
\end{tikzpicture}
\caption{Reduction of oscillation at $x_0$} 
\label{fig:x0}
\end{figure}

It is sufficient to 
show that 
  \[
   |u(x)-u(x_0)|\leq c\tfrac{L(|x-x_0|)^{-\beta}}{L(r)^{-\beta}}\quad \text{ for any }\ x\in B_r\,.
  \]
  
  Define  $r_n=L^{-1}(a^{n-1}L(r/2))$ for $n\in \N$, 
  where $a>2$ will be  chosen in the course of the proof independently of $r$, $u$ and 
$f$. 
  We will construct a nondecreasing sequence $(c_n)_{n\in \N}$ and 
non-increasing sequence 
$(d_n)_{n\in \N}$ of positive numbers so that 
  \begin{equation}\label{eq:todo}
   c_n\leq u(x)\leq d_n \text{ for all }\ x\in B_{r_n}(x_0)\quad \text{ and }\quad d_n-c_n\leq b^{-n+1}\,,
  \end{equation}
where $b=\frac{2}{2-\theta} \in (1,2) $ and $\theta \in (0,1)$ will be chosen later 
independently of $r, u$ and 
$f$.
   This will be enough, since for $r_{n+1}\leq |x-x_0|< r_n$ we will then have 
   \begin{align*}
    |u(x)-u(x_0)|&\leq b^{-n+1}=b\left(\tfrac{1}{a^{n}}\right)^{\frac{\ln{b}}{\ln{a}}}=b\left(\tfrac{L(r/2)}{a^{n}L(r/2)}\right)^{\frac{\ln{b}}{\ln{a}}}\\
      &=b\left(\tfrac{L(r/2)}{L(r_{n+1})}\right)^{\frac{\ln{b}}{\ln{a}}}\leq b\left(\tfrac{L(r/2)}{L(|x-x_0|)}\right)^{\frac{\ln{b}}{\ln{a}}}\\
      &\leq b \big(\tfrac{2^ 
\gamma}{c_{L}}\big)^{\frac{\ln{b}}{\ln{a}}} \left(\tfrac{L(|x-x_0|)}{L(r)}\right)^{- 
\frac{\ln{b}}{\ln{a }}},
   \end{align*}
   where in the last inequality we have used \autoref{lem:L}.

   We prove (\ref{eq:todo}) and construct sequences $(c_n)$ and $(d_n)$ 
inductively.
   We set    
   \[
         c_1:=\inf_{\R^d}u\qquad \text{ and }\qquad d_1:=c_1+1 \,.     
   \]
  Let $n\in \N$, $n \geq 2$. Assume that $c_k$ and $d_k$ have been constructed for $k\leq 
n$ and that 
(\ref{eq:todo}) holds for 
$k\in \N, k\leq n$. We are now going to construct $c_{n+1}$ and 
$d_{n+1}$.
  
  Set $m:=\frac{c_n+d_n}{2}$. By (\ref{eq:todo}) it follows for $x \in B_{r_n}$ 
  \[
   u(x)-m\leq \tfrac{1}{2}(d_n-c_n)
   \leq\tfrac{1}{2}b^{-n+1}\,.
  \]
  Define a function $v\colon\R^d\rightarrow \R$ by $v(x):=2 
b^{n-1}(u(x_0+x)-m)$. Then 
$v(x)\leq 1$ for $x\in B_{r_n}$ 
and $Av=2 b^{n-1} f$ in $B_{r_n}$.
  
  Assume that $\mu(\{x\in B_{r_n}\setminus B_{r_{n+1}}\colon v(x)\leq 0\})\geq 
\frac{1}{2}\mu(B_{r_n}\setminus B_{r_{n+1}})$. We recall that the ball $B_{r/2}$ has 
center $0$ and the balls $B_{r_n}$ have center $x_0$. For $x\in B_{r/2}\setminus 
B_{r_n}$ there exists $k\in \N$, $k\leq n-1$ such that $r_{k+1}\leq |x|<r_k$. Then by 
$(\ref{eq:todo})$ we have 
  \begin{align*}
   v(x)&=2b^{n-1}(u(x_0+x)-m)\leq 2b^{n-1}(d_k-m)\leq 2b^{n-1}(d_k-c_k)\\
   &= 2b^{n-k}=2b\left(\tfrac{a^{n-1}L(r/2)}{a^{k}L(r/2)}\right)^{\frac{\ln{b}}{\ln{a}}}=2b\left(\tfrac{L(r_n)}{L(r_{k+1})}\right)^{\frac{\ln{b}}{\ln{a}}}\\
   &\leq 2b\left(\tfrac{L(r_n)}{L(|x|)}\right)^{\frac{\ln{b}}{\ln{a}}}\,.
  \end{align*}
  If $x\in B_{r/2}^c$, then $L(|x|)\leq L(r/2)$ and  $u(x_0+x)-m\leq d_1-c_1=1$; hence 
  \[
  	v(x)\leq 2 b^{n-1}
= 2 \left(\tfrac{a^{n-1}L(r/2)}{L(r/2)}\right)^{\frac{\ln{b}}{\ln{a}}} 
= 2 \left(\tfrac{L(r_n)}{L(r/2)}\right)^{\frac{\ln{b}}{\ln{a}}} 
\leq 2 \left(\tfrac{L(r_n)}{L(|x|) }\right)^{\frac{\ln{b}}{\ln{a}}}\,.
  \]
  
We want to apply \autoref{lem:expans1} with $r=r_{n+1}$. Note that 
$r_n=\varphi_a(r_{n+1})$. In order to 
apply \autoref{lem:expans1} we need to verify that $2 b^{n-1} |f| \leq 
L(\varphi(r_{n+1})) 
= a^{n-1} L(r/2)$. But this holds true because $|f| \leq \tfrac12  L(r/2)$ and 
$(b/a)^{n-1} \leq b/a \leq 2 a^{-1} \leq 1$. Thus we obtain that for some $a>2$ and 
$\theta\in (0,1)$, not depending on $v$ and $r$, $v(x)\leq 1-\theta$ on 
$B_{r_{n+1}}$. Going back to $u$, we deduce
  \[
  	u(x)\leq \frac{1-\theta}{2} \, b^{1-n}+\frac{c_n+d_n}{2}\qquad \text{ for }\ 
x\in B_{r_{n+1}}(x_0)\,.
  \]
We take $c_{n+1}=c_n$ and 
$d_{n+1}=\min\{d_n, \frac{1-\theta}{2} \, b^{1-n}+\frac{c_n+d_n}{2} \}$. This choice 
implies $d_{n+1}-c_{n+1} \leq b^{-n}$. 

In the case $\mu(\{x\in B_{r_n}\setminus B_{r_{n+1}}\colon v(x)\leq 0\})< 
\frac{1}{2}\mu(B_{r_n}\setminus B_{r_{n+1}})$ we repeat the previous argument with $-v$ 
instead of $v$ and deduce  
   \[
  	u(x)\geq - \frac{1-\theta}{2} \, 
b^{1-n}+\frac{c_n+d_n}{2} \qquad \text{ for 
}\quad x\in B_{r_{n+1}}(x_0)\,.
  \]
This time we choose $c_{n+1}= \max\{c_n, - \frac{1-\theta}{2} \, 
b^{1-n}+\frac{c_n+d_n}{2}\}$ and $d_{n+1}=d_n$. Finally, we set 
$\beta:=\frac{\ln{b}}{\ln{a}}$.

We have completed the proof in the special case where conditions \eqref{eq:A1} 
and \eqref{eq:K-Levy-bound} hold true with $R_0 = + 
\infty$. Now, let us  assume that 
conditions \eqref{eq:A1} and \eqref{eq:K-Levy-bound} hold true for some 
$R_0 < + \infty$, i.e., there is $\ell\colon 
(0,R_0)\rightarrow 
(0,+\infty)$ satisfying $\int_0^{R_0} \ell(s)\frac{ds}{s} = + 
\infty$ and conditions \eqref{eq:ell1}, \eqref{eq:ell2} and 
\eqref{eq:ell3} for some $c_L \in (0,1)$, $c_U \geq 1$, 
and $\gamma\in (0,2)$. We define $\widetilde{\ell}\colon 
(0,+\infty)\rightarrow (0,+\infty)$ by
\begin{align*}
 \widetilde{\ell}(s) = 
\begin{cases}
\ell(s) \quad &\text{ for } 0 < s < R_0 \,, \\
(s-\frac{R_0}{2})^{-\gamma} &\text{ for } s \geq R_0 \,.
\end{cases}
\end{align*}

\begin{lemma}
The function $\widetilde{\ell}$ satisfies \eqref{eq:ell1}, \eqref{eq:ell2} and 
\eqref{eq:ell3} with $R_0$ being replaced by $+ \infty$ and for some 
$\widetilde{c_L} \in (0,1)$, 
$\widetilde{c_U} \geq 1$, $\gamma\in (0,2)$.
\end{lemma}

\proof
Assume $R_0 > 0$. Condition \eqref{eq:ell1} obviously holds true. When checking 
\eqref{eq:ell2} 
and \eqref{eq:ell3} we need to consider several cases. Both conditions hold 
true for $\lambda r < R_0$ because the functions $\widetilde{\ell}$ and $\ell$ 
coincide in this range. For $\lambda r \geq R_0, r \geq R_0$ the conditions can 
easily be verified. The only challenging case is $\lambda r \geq R_0, r < R_0$. 
Let us verify \eqref{eq:ell2} in this case, i.e., show that 
\begin{align}\label{eq:ell2tildecheck}
\widetilde{\ell}(r\lambda) \geq \widetilde{c_L} \lambda^{-\gamma} 
\widetilde{\ell}(r) \quad \Leftrightarrow \quad (r\lambda - 
\tfrac{R_0}{2})^{-\gamma} 
\geq \widetilde{c_L} \lambda^{-\gamma} \ell (r) \,,
\end{align}
which is equivalent to
\begin{align*}
\frac{(r \lambda - \frac{R_0}{2})^{- \gamma}} {(r \lambda)^{-\gamma}} 
\, r^{-\gamma} \geq  \widetilde{c_L} \ell (r) \quad \text{ and } 
\quad \Big(\frac{r \lambda}{r \lambda - 
\frac{R_0}{2}}\Big)^\gamma \, r^{-\gamma} \geq  \widetilde{c_L} \ell (r)\,.
\end{align*}
Since the fraction on the left-hand side is bounded in $[1,2]$ for 
$r\lambda \geq R_0$, it is sufficient to show $r^{-\gamma} >  \widetilde{c_L} 
\ell (r)$ in the case $r < R_0 \leq \lambda r$ for some $\widetilde{c_L} 
\in (0,1)$. Let us prove this assertion for two cases separately. In the case 
$\frac{R_0}{2} < r < R_0 \leq \lambda r$ we conclude:
\begin{align*}
\ell(r) &= \ell(\tfrac{R_0}{2}) \tfrac{\ell(r)}{\ell(\frac{R_0}{2})} 
\overset{\eqref{eq:ell3}}{\leq} \ell(\tfrac{R_0}{2}) c_U 
2^dR_0^{-d}r^d \leq \ell(\tfrac{R_0}{2}) c_U 2^d 
\leq  \ell(\tfrac{R_0}{2}) c_U 2^d \tfrac{r^{-\gamma}}{R_0^{-\gamma}}
\end{align*}
which proves the assertion. In the case $r < \frac{R_0}{2} < R_0 \leq \lambda 
r$ we proceed as follows:
\begin{align*}
\ell(r) &= \tfrac{\ell(\frac{R_0}{2})}{\frac{\ell(\frac{R_0}{2})}{\ell(r)}} 
\overset{\eqref{eq:ell2}}{\leq} \tfrac{\ell(\frac{R_0}{2})}{c_L 
(\frac{R_0}{2r})^{-\gamma}} = c_L (\tfrac{R_0}{2})^{\gamma} 
\ell(\tfrac{R_0}{2}) r^{-\gamma} \,. 
\end{align*}
Thus we have shown 
\[ \widetilde{\ell}(r\lambda) \geq \widetilde{c_L} \lambda^{-\gamma} 
\widetilde{\ell}(r) \quad \text{ for all } r>0 \text{ and }  \lambda > 1 \]
for an appropriate choice of $\widetilde{c_L}$. In other words, the function 
$\widetilde{\ell}$ satisfies condition \eqref{eq:ell2} with $R_0$ being 
replaced by $+ \infty$. It remains to show 
\[ \widetilde{\ell}(r\lambda) \leq \widetilde{c_U} \lambda^{d} 
\widetilde{\ell}(r) \quad \text{ for all } r>0 \text{ and }  \lambda > 1 \]
for an appropriate choice of $\widetilde{c_U}$. As explained above, it remains 
to show this estimate in the case $r < R_0 
\leq \lambda r$.
In the case $r\leq \frac{R_0}{2}<R_0\leq \lambda r$ we obtain
\begin{align*}
\frac{\ell(\lambda r)}{\ell(r)}&=\frac{\left(\lambda 
r-\frac{R_0}{2}\right)^{-\gamma}}{\ell(\frac{R_0}{2})}\frac{\ell(\frac{R_0}{2})}
{\ell(r)}\overset{\eqref{eq:ell3}} 
{\leq} 2^{\gamma} R_0^{-\gamma} 
\ell(\tfrac{R_0}{2})^{-1}c_U\left(\tfrac{R_0}{2r}\right)^d\\
&\leq 2^{\gamma} R_0^{-\gamma} \ell(\tfrac{R_0}{2})^{-1}c_U \lambda^d\,.
\end{align*}
If $\frac{R_0}{2}<r<R_0\leq \lambda r$, 
\begin{align*}
\frac{\ell(\lambda r)}{\ell(r)}&=\frac{\left(\lambda 
r-\frac{R_0}{2}\right)^{-\gamma}}{\ell(\frac{R_0}{2})}\frac{\ell(\frac{R_0}{2})}
{\ell(r)}\overset{\eqref{eq:ell2}} 
{\leq} 2^{\gamma} 
R_0^{-\gamma} \ell(\tfrac{R_0}{2})^{-1}c_L^{-1}\left(\tfrac{R_0}{2r}\right)^{ 
\gamma}\\
&\leq 2^{\gamma} R_0^{-\gamma} \ell(\tfrac{R_0}{2})^{-1}c_L^{-1}\leq 
2^{\gamma} R_0^{-\gamma} \ell(\tfrac{R_0}{2})^{-1}c_L^{-1} \lambda^{d}\,.
\end{align*}
\qed

Next, define a modified kernel function $\widetilde{K}(x,h)$ by 
\begin{align*}
 \widetilde{K}(x,h) = 
\begin{cases}
K(x,h) \quad &\text{ for } 0 < |h| < R_0 \,, \\
|h|^{-d} (|h|-\tfrac{R_0}{2})^{-\gamma} &\text{ for } |h| \geq R_0 \,,
\end{cases}
\end{align*}
with $\gamma$ as in \eqref{eq:ell2}. Let us denote the integrodifferential 
operator corresponding to $\widetilde{K}$ 
by $\widetilde{A}$. Since $u$ solves $A u = f$, the function $u$ also solves
\[ \widetilde{A} u = f + (\widetilde{A} - A ) u := \widetilde{f} \,. \]
Due to the definition of $\widetilde{K}$, the image of $u$ under the operator 
$(\widetilde{A} - A )$ is a bounded function, hence $\widetilde{f}$ is 
a bounded function. Now we apply the previous proof ($R_0 = + \infty$) and 
obtain for every $r >0$
\begin{align}\label{eq:theo-estim-R0infty}
  \sup_{x,y\in B_{r/4}}\tfrac{|u(x)-u(y)|}{L(|x-y|)^{-\beta}}\leq 
c L(r)^{\beta}\|u\|_\infty + c L(r)^{\beta-1} \|f\|_{L^\infty(B_r)} + c 
L(r)^{\beta-1} \|(\widetilde{A} - A ) u\|_{L^\infty(B_r)} 
\end{align}
for some positive constant $c>0$. Note that
\[ \|(\widetilde{A} - A ) u\|_{L^\infty(B_r)} \leq 2  \|u\|_\infty \sup_{x \in 
B_r} \int\limits_{\R^d\setminus B_{R_0}} \left(K(x,h) + \widetilde{K}(x,h)\right) \d h \leq 
c' \|u\|_\infty
\]
for some constant $c'>0$. Since $L(\frac{R_0}{2}) \leq L(r)$ for $0 < r \leq 
\frac{R_0}{2}$ we conclude
\begin{align*}
\|(\widetilde{A} - A ) u\|_{L^\infty(B_r)} \leq c' L(\tfrac{R_0}{2})^{-1} L(r)
\|u\|_\infty
\end{align*}
and finally
\begin{align*}
  \sup_{x,y\in B_{r/4}}\tfrac{|u(x)-u(y)|}{L(|x-y|)^{-\beta}}\leq 
(c+c' L(\tfrac{R_0}{2})^{-1}) L(r)^{\beta}\|u\|_\infty + c L(r)^{\beta-1} 
\|f\|_{L^\infty(B_r)}\,. 
\end{align*}
\qed


\section{An approach to regularity via stochastic processes}\label{sec:intro-prob}

As explained in the introduction, the aim of this section is to provide an 
alternative approach to \autoref{theo:main-analysis} using stochastic 
processes. First, let us formulate our assumptions and results. As in the first 
part, we assume $0 \leq 
\alpha < 2$ and $K\colon \R^d\times(\R^d\setminus\{0\})\rightarrow [0,+\infty)$ 
to be a
measurable function satisfying \eqref{eq:K-Levy-bound} and the symmetry condition 
$K(x,h)=K(x,-h)$ for all $x,h$. Instead of condition \eqref{eq:A1} we assume the 
following condition:
\begin{align}
\label{eq:K3} \tag{K} & \ \ \  \kappa^{-1} \, \frac{\ell(|h|)}{|h|^d}\leq
K(x,h)\leq
\kappa \, \frac{\ell(|h|)}{|h|^d} \quad \text{ for } 0 < |h|\leq 1 \,,
\end{align}
where $\kappa>1$ and $\ell\colon
(0,1)\rightarrow (0,+\infty)$ is locally bounded and varies
regularly at zero with index $-\alpha\in (-2,0]$. Possible examples could be
$\ell(s)=1$,
$\ell(s)=s^{-3/2}$ and $\ell(s)=s^{-\beta} \ln(\tfrac{2}{s})^{2}$ for some
$\beta \in (0,2)$, see \autoref{sec:appendix} for a more detailed discussion.  

These assumptions on $K(x,h)$ differ slightly from the ones in 
\autoref{sec:intro-main}. Concerning the behavior of $K(x,h)$ for small values 
of $|h|$, these assumptions are slightly more restrictive. We suppose that 
there exists a strong Markov process $X=(X_t,\P_x)$ with
trajectories that are right continuous with left limits associated with
${\opA}$ in the sense that for every $x \in \R^d$
\begin{itemize}
\item[(i)] $\P_x(X_0=x)=1$;
\item[(ii)] for any $f\in C_b^2(\R^d)$  the process $
\big(f(X_t)-f(X_0)-\int_0^t
{\opA}f(X_s) \d s \big)_{t \geq 0}$ is a martingale under $\P_x$.
\end{itemize}

Note that the existence of such a Markov process comes for free in the
case when $K(x,h)$ is independent of $x$, see \autoref{sec:intro-prob}. In
the general case it has been established by many authors in different general
contexts, see the discussion in \cite{AbKa09}. Denote by $\tau_A=\inf\{t>0|\,
X_t\not\in A\}$,
$T_A=\inf\{t>0|\, X_t \in A\}$  the first exit time resp. hitting time
of the
process $X$ for a measurable set $A\subset \R^d$.

\begin{definition}
 A bounded function $u\colon \R^d\rightarrow \R$ is said to be harmonic in an
open subset $D\subset \R^d$ with respect to $X$ (and ${\opA}$) if for
any bounded open set  $B\subset \overline{B}\subset D$ the stochastic process
$(u(X_{\tau_B\wedge t}))_{t\geq 0}$ is a $\P_x$-martingale for every $x\in
\R^d$.
\end{definition}

Before we can formulate our results we need to introduce an additional quantity.
Note
that \eqref{eq:K-Levy-bound} and \eqref{eq:K3} imply that $\int_0^1 s \,
\ell(s) \, \d s \leq c$ holds for some constant $c>0$. Let
$L\colon (0,1)\rightarrow (0,+\infty)$ be defined by
$ L(r)=\int_r^1\frac{\ell(s)}{s}\, \d s$. The function $L$ is well
defined because $L(r) \leq r^{-2} \int_r^1 s^2 \frac{\ell(s)}{s}\, \d s \leq c
r^{-2}$.  See the table in \autoref{sec:intro-main} for several examples.  We 
note that
the function $L$ is always decreasing. 

{\bf Remark:} The definition of $L$ here is different from the definition in 
\autoref{sec:intro-main}. The reason is that, here we are able to work without
specific assumptions on $K(x,h)$ for large values of $|h|$.

The analogous result to \autoref{theo:main-analysis}, which we prove with 
probabilistic techniques, is the following theorem.

\begin{theorem}\label{theo:main-prob} There exist constants $c>0$ and $\gamma\in 
(0,1)$
so that for all $r\in (0,\frac{1}{2})$ and $x_0\in \R^d$
\begin{align}\label{eq:main_reg}
|u(x)-u(y)|\leq c\|u\|_\infty
\frac{L(|x-y|)^{-\gamma}}{L(r)^{-\gamma}},\ \ x,y\in
B_{r/4}(x_0)
\end{align}
for all bounded functions $u\colon \R^d\rightarrow \R$ that are harmonic in
$B_r(x_0)$ with respect to $\opA$.
\end{theorem}

Let us comment on this result. It is important to note that the result
trivially holds if the function $L:(0,1) \to (0,+\infty)$ satisfies
$\lim\limits_{r\to 0+}L(r) < +\infty$. This is equivalent to the condition
\begin{align} \label{eq:levy-finite}
\int\limits_{B_1} \frac{\ell(|h|)}{|h|^d} \, \d h < +\infty \,,
\end{align}
which, in the case $K(x,h) = k(h)$, means that the L\'{e}vy measure is finite. 
Thus, for the proof, we can concentrate on cases where \eqref{eq:levy-finite}
does not hold. One could say that our result holds true up to
and across the phase boundary determined by whether the kernel
$K(x,\cdot)$ is integrable (finite L\'{e}vy measure) or not. 

Furthermore, note that the main result of \cite{BaLe02} is implied by
\autoref{theo:main-prob} since
the choice $\ell(s)=s^{-\alpha}$, $\alpha \in (0,2)$, leads to
$L(r)\asymp r^{-\alpha}$. Given the whole spectrum of possible operators covered
by
our approach, this choice is a very specific one. It allows to use scaling
methods in the usual way which are not at our disposal here.
\autoref{tab:choices_l} in \autoref{sec:intro-main} contains several admissible
examples. Note that \eqref{eq:main_reg}
becomes trivial if $L(0) < +\infty$.  

The main ingredient in the proof of \autoref{theo:main-prob} is a new version of
\autoref{prop:hitting-prop} which we provide now. For $r\in (0,1)$ we
define a measure $\mu$ by 
\begin{equation}\label{eq:measure_mu_null}
 \mu (dx)=\frac{\ell(|x|)}{L(|x|) |x|^{d}} \, \mathbbm{1}_{B_1}(x) \d x\,.
\end{equation}
Moreover, for $a>1$, we define a function $\varphi_a:(0,1) \to (0,1)$ by
$\varphi_a(r)=L^{-1}(a^{-1} L(r))$. The following result is our
modification of \autoref{prop:hitting-prop}. 

\begin{proposition}\label{prop:hitting_new}
There exists a constant $c>0$ such that for all $a>1$, $r\in
(0,\frac{1}{2})$ and measurable sets $A\subset B_{\varphi_a(r)}\!\setminus\!
B_r$ with $\mu (A)\geq
\frac{1}{2}\mu (B_{\varphi_a(r)}\!\setminus\! B_r)$
 \[
  \P_x(T_A<\tau_{B_{\varphi_a(r)}}) \geq \P_x(X_{\tau_{B_r}}\in A)\geq
c \, \tfrac{\ln{a}}{a}
 \]
holds true for all $x\in B_{r/2}$.
\end{proposition}

The novelty of \autoref{prop:hitting_new} is the definition and use of
the measure $\mu$. It allows us to deal with the classical cases as well as
with critical cases, e.g. given by $K(x,h)\asymp |h|^{-d}
\mathbbm{1}_{B_1}(h)$. 

Note that we use the notation $f(r)\asymp g(r)$ to 
denote that the ration $f(r)/g(r)$ stays between two positive constants as $r$
converges to some value of interest.

%

\section{Probabilistic estimates}\label{sec:prob_estim}
\begin{proposition}\label{prop:estp} There exists a constant $C_1>0$ such that
for $x_0\in \R^d$, $r\in
(0,1)$ and $t>0$ 
\[
\P_{x_0}(\tau_{B_r(x_0)}\leq t)\leq C_1t\,L(r) \,.
\]
\end{proposition}
\proof
Let $x_0\in \R^d$, $0<r<1$ and let $f\in C^2(\R^d)$ be a positive function such 
that
\[
f(x)=\left\{\begin{array}{cl}
|x-x_0|^2, & |x-x_0|\leq \frac{r}{2}\\
r^2,& |x-x_0|\geq r
\end{array}\right.
\]
and for some $c_1 >0$ 
\[
|f(x)|\leq c_1 r^2, \ \ \left|\frac{\partial f}{\partial x_i}(x)\right|\leq c_1 
r \ \ \textrm{ and }
\ \ \left|\frac{\partial^2 f}{\partial x_i \partial x_j}(x)\right|\leq c_1.
\]
By the optional stopping theorem we get
\begin{align}\label{eq:optional_stopping}
\E_x f(X_{t\wedge \tau_{B_r(x_0)}})-f(x_0)=\E^x \int_0^{t\wedge
\tau_{B_r(x_0)}}\mathcal{L}f(X_s)\,
ds,\ \ t>0.
\end{align}
Let $x\in B_r(x_0)$. We estimate $\mathcal{L}f(x)$ by splitting the integral in
\eqref{eq:def_A} into three parts.
\begin{align*}
\int_{B_r} &(f(x+h)-f(x)-\nabla f(x)\cdot h \mathbbm{1}_{\{|h|\leq
1\}})K(x,h)\d h\\
& \leq c_2\int_{B_r}|h|^2K(x,h)\d h\leq c_2\kappa
\int_{B_r}|h|^{2-d}\ell(|h|)\d h \leq c_3 r^{2}\ell(r),
\end{align*}
where in the last line we have used Karamata's theorem, see property
\eqref{rem:karamata} in \autoref{sec:appendix}. On the other hand, on $B_r^c$ we
have
\begin{align*}
\int_{B_r^c} &(f(x+h)-f(x))K(x,h)\d h\leq 2\|f\|_\infty
\int_{B_r^c}K(x,h)\d h\\
& \leq 2\|f\|_\infty \left(\kappa\int_{B_1\setminus
B_r}|h|^{-d}\ell(|h|)\d h+\int_{B_1^c}K(x,h)\d h\right) \leq c_4r^{2}L(r)\d r
\,,
\end{align*}
where we applied property \eqref{rem:L_notzero} from \autoref{sec:appendix}.
Last, we estimate
\begin{align*}
\left|\int_{B_1\setminus B_r}h\cdot \nabla f(x) K(x,h)\d h\right|&\leq
c_1 r\int_{B_1\setminus B_r} |h| K(x,h)\d h\\
&\leq c_1 \kappa r \int _
{B_1\setminus B_r} |h|^{-d+1}\ell(|h|)\d h \leq c_5 r^2\ell(r),
\end{align*}
by Karamata's theorem again. Therefore, by property \eqref{rem:l_versus_L} from
\autoref{sec:appendix} we conclude that
there is a constant $c_6>0$ such that for all $x\in B_r(x_0)$ and $r\in (0,1)$
we have
\begin{equation}\label{eq:prop1e1}
\mathcal{L}f(x)\leq c_6 r^2 L(r).
\end{equation}
Let us look again at \eqref{eq:optional_stopping}. On
$\{\tau_{B_r(x_0)}\leq t\}$ we have $X_{t\wedge \tau_{B_r(x_0)}}\in
B_r(x_0)^c$ and so
$f(X_{t\wedge \tau_{B_r(x_0)}})\geq r^2$. Thus, by  (\ref{eq:prop1e1}) and
\eqref{eq:optional_stopping} we get
\[
\P_{x_0}(\tau_{B_r(x_0)}\leq t)\leq c_6 L(r)t.
\]
\qed

\begin{proposition}\label{prop:prob-2}
	There are constants $C_2>0$ and $C_3>0$ such that for $x_0\in \R^d$
	\[
		\sup_{x\in \R^d}\E_x\tau_{B_r(x_0)}\leq
\frac{C_2}{L(r)}\,,\quad r\in
(0,1/2)
	\]
and
	\[
		\inf_{x\in B_{r/2}(x_0)}\E_x\tau_{B_r(x_0)}\geq
\frac{C_3}{L(r)}\,,\quad r\in
(0,1)
	\]

\end{proposition}
\proof
	The proof is similar to the proof of the exit time estimates in 
\cite{BaLe02}. 

      (a) First we prove the upper estimate for the exit time.
	Let $x\in\R^d$, $r\in (0,1/2)$ and let
	\[
		S=\inf\{t>0|\, |X_t-X_{t-}|>2r\}
	\]
	be the first time of a jump larger than $2r$. With the help of the
L\' evy system formula (see \cite[Proposition 2.3]{BaLe02}) and 
\eqref{eq:K3} we
can deduce
{\allowdisplaybreaks 
	\begin{align}
		\P_x(S\leq L(r)^{-1}) & =\E_x\sum_{t\leq L(r)^{-1}\wedge S}
\mathbbm{1}_{\{|X_t-X_{t-}|>2r\}} =\E_x\int\limits_0^{L(r)^{-1}\wedge
S}\int\limits_{B_{2r}^c}K(X_s,h)\d h\d s \nonumber\\
		&\geq c_1\E_x[L(r)^{-1}\wedge 
S]\int\limits_{2r}^1\frac{\ell(t)}{t}\d t\,.\label{eq:prob-6}
	\end{align}
	}
	Since $L$ is regularly varying at zero, 
	\begin{align*}
		\E_x[L(r)^{-1}\wedge S]& \geq L(r)^{-1}\P_x(S>L(r)^{-1}) \geq
c_2L(2r)^{-1}\big(1-\P_x(S\leq L(r)^{-1})\big)\,
	\end{align*}
	and so it follows from (\ref{eq:prob-6}) that
	\begin{equation}\label{eq:prob-7}
		\P_x(S\leq L(r)^{-1})\geq c_3
	\end{equation}
	with $c_3=\frac{c_1c_2}{c_1c_2+1} \in (0,1)$. The strong Markov
property and (\ref{eq:prob-6}) lead to
	\[
		\P_x(S>m L(r)^{-1})\leq (1-c_3)^m,\ \ m\in \N\,.
	\]
	Since $\tau_{B_r(x_0)}\leq S$, 
	\begin{align*}
		\E_x\tau_{B_r(x_0)}\leq \E_x S &\leq
L(r)^{-1}\sum_{m=0}^\infty (m+1)\P_x(S>L(r)^{-1}m)\\
		&\leq L(r)^{-1}\sum_{m=0}^\infty (m+1)(1-c_3)^m\,.
	\end{align*}

(b) Now we prove the lower estimate of the exit time. Let $r\in
(0,1)$ and
$y\in B_{r/2}(x_0)$.
By \autoref{prop:estp},
\[
  \P_y(\tau_{B_r(x_0)}\leq t)\leq \P_y(\tau_{B_{r/2}(y)}\leq t)\leq
C_1tL(r/2),\quad t>0\,,
\]
since $B_{r/2}(y)\subset B_r(x_0)$. Choose $t=\frac{1}{2C_1L(r/2)}$. Then
\begin{align*}
 \E_y\tau_{B_r(x_0)}&\geq \E_y[\tau_{B_r(x_0)};\tau_{B_r(x_0)}>t]\geq
t\P_y(\tau_{B_r(x_0)}>t)\\
		    &\geq t(1-C_1L(r/2)t)=\frac{1}{4C_1L(r/2)}\,.
\end{align*}
By \eqref{rem:l-reg_L-reg} from \autoref{sec:appendix} we know that
$L$ is regularly varying at zero. Hence there is a
constant $c_1>0$ such that $L(r/2)\leq c_1 L(r)$ for all $r\in (0,1/2)$.
Therefore
$\E_y\tau_{B_r(x_0)}\geq \frac{1}{4C_1c_1L(r)}$.

\qed

\begin{proposition}\label{prop:prob-3}
	There is a constant $C_4>0$ such that for all $x_0\in \R^d$ and $r,s\in
(0,1)$ satisfying $2r<s$
	\[
		\sup_{x\in B_r(x_0)}\P_x(X_{\tau_{B_r(x_0)}}\not\in
B_s(x_0))\leq C_4 \frac{L(s)}{L(r)}\,.
	\]
\end{proposition}
\proof
	Let $x_0 \in \R^d$, $r,s \in (0,1)$ and $x\in B_r(x_0)$. Set
$B_r:=B_r(x_0)$. By the L\' evy system formula, for $t>0$
	\begin{align*}
		\P_x(X_{\tau_{B_r}\wedge t}\not\in B_s)&=\E_x\sum\limits_{v\leq
\tau_{B_r}\wedge t} \mathbbm{1}_{\{X_{v-}\in B_r,X_v\in
B_s^c\}} =\E_x\int\limits_0^{\tau_{B_r}\wedge
t}\int\limits_{B_s^c} K(X_v,z-X_v)\d z\d v \,.
	\end{align*}
	
	Let  $y\in B_r$. Since $s\geq 2r$, it follows that
$B_{s/2}(y)\subset B_s$ and hence
	\begin{align*}
		\int\limits_{B_s^c} K(y,z-y)\d z&\leq
\int\limits_{B_{s/2}(y)^c} K(y,z-y)\d z \leq
c_1\int_{s/2}^1 \frac{\ell(u)}{u}\d u+c_2 \leq c_3L(s)\,.
	\end{align*}
	where in the last inequality we have used that $L$ varies
regularly at zero and that $\lim\limits_{r\to 0+}L(r)> 0$, see
\eqref{rem:L_notzero} in \autoref{sec:appendix}. 
	
	The above considerations together with \autoref{prop:prob-2} imply
	\[
		\P_x(X_{\tau_{B_r}\wedge t}\not\in B_s)\leq c_3 L(s)\E_x 
\tau_{B_r}\leq c_4 \frac{L(s)}{L(r)}\,.
	\]
	Letting $t\to\infty$ we obtain the desired estimate.
\qed

For $x_0\in \R^d$ and $r\in (0,1)$ we define the following measure
\begin{equation}\label{eq:measure_mu}
 \mu_{x_0}(dx)=\frac{\ell(|x-x_0|)}{L(|x-x_0|)}\,|x-x_0|^{-d}
\mathbbm{1}_{\{|x-x_0|<1\}}\d x\,.
\end{equation}

Define $\varphi_a(r)=L^{-1}(\frac{1}{a}L(r))$ for $r\in (0,1)$ and $a>1$. The
following property is important for the construction below:
\begin{align} 
r = L^{-1}(L(r)) \leq L^{-1}(\tfrac{1}{a}L(r)) = \varphi_a(r)
\,. 
\end{align}

Now we can prove a Krylov-Safonov type hitting estimate which
includes \autoref{prop:hitting_new} as a special case. 

\begin{proposition}\label{prop:hitting_general}
There exists a constant $C_5>0$ such that for all $x_0\in \R^d$, $a>1$, $r\in
(0,\frac{1}{2})$ and
$A\subset B_{\varphi_a(r)}(x_0)\setminus B_r(x_0)$ satisfying
$\mu_{x_0}(A)\geq
\frac{1}{2}\mu_{x_0}(B_{\varphi_a(r)}(x_0) \setminus
B_r(x_0))$
 \[
  \P_y(T_A<\tau_{B_{\varphi_a(r)}(x_0)})\geq \P_y(X_{\tau_{B_r(x_0)}}\in A)\geq
C_5\frac{\ln{a}}{a}\,,\quad y\in B_{r/2}(x_0)\,.
 \]
\end{proposition}

\proof
Consider $x_0\in \R^d$, $a>1$, $r\in
(0,\frac{1}{2})$ and a set 
$A\subset B_{\varphi_a(r)}(x_0)\setminus B_r(x_0)$ satisfying
$\mu_{x_0}(A)\geq
\frac{1}{2}\mu_{x_0}(B_{\varphi_a(r)}(x_0) \setminus
B_r(x_0))$. Set $\mu:=\mu_{x_0}$, $\varphi:=\varphi_a$,  $B_s:=B_s(x_0)$ and
let $y\in
B_{r/2}$. The first
inequality follows from 
$\{X_{\tau_{B_r}}\in A\}\subset \{T_A<\tau_{B_{\varphi(r)}}\}$ since $A\subset
B_{\varphi(r)}\setminus
B_r$.

By the L\' evy system formula, for $t>0$, 
	\begin{align}\label{eq:ks-0}
		\P_y(X_{\tau_{B_r}\wedge t}\in A)&=\E_y\sum\limits_{s\leq 
\tau_{B_r}\wedge
t} \mathbbm{1}_{\{X_{s-}\in B_r,X_s\in A\}}	
=\E_y\int\limits_0^{\tau_{B_r}\wedge t}\int\limits_{A} K(X_s,z-X_s)\d z\d s\,.
	\end{align}

Since $|z-x|\leq |z-x_0|+|x_0-x|\leq |z-x_0|+r\leq 2|z-x_0|$ for $x\in B_r$ and 
$z\in B_r^c$,
\begin{equation}\label{eq:ks-1}
 \E_y\int\limits_0^{\tau_{B_r}\wedge t}\int\limits_{A} K(X_s,z-X_s)\d z\d s\geq
c_1\E_y[\tau_{B_r}\wedge t]\int_A \frac{\ell(|z-x_0|)}{|z-x_0|^d}\d z\,,
\end{equation}
where we have used property \eqref{rem:potter} given in
\autoref{sec:appendix}.

Since $L$ is decreasing,
\begin{align}
 \int_A \frac{\ell(|z-x_0|)}{|z-x_0|^d}\d z&=\int_A
L(|z-x_0|)\mu(dz) \geq
L(\varphi(r))\mu(A)\geq
\frac{L(r)}{2a}\mu(B_{\varphi(r)}\setminus B_r)\,.\label{eq:ks-2}
\end{align}
Noting that 
\[\mu(B_{\varphi(r)}\setminus B_r)=c_2\int_r^{\varphi(r)}
\frac{1}{L(s)}\frac{\ell(s)\d s}{s}=-c_2\ln L(s)|_r^{\varphi(r)}=c_2\ln a\,,\]
we conclude from
(\ref{eq:ks-0})--(\ref{eq:ks-2}) that
\[
  \P_y(T_A<\tau_{B_{\varphi_a(r)}(x_0)})\geq c_3
L(r)\frac{\ln{a}}{a}\E_y[\tau_{B_r}\wedge t]\,.
\]
Letting $t\to\infty $ and using the lower bound in \autoref{prop:prob-2} we get
\begin{align*}
 \P_y(T_A<\tau_{B_{\varphi_a(r)}(x_0)}) \geq c_3
L(r)\,\frac{\ln{a}}{a}\,\E_y\tau_{B_r} \geq
c_3L(r)\,\frac{\ln{a}}{a}\,C_3L(r)^{-1}=c_3C_3\frac{\ln{a}}{a}\,.
\end{align*}
\qed

\section{Proof of \autoref{theo:main-prob}}\label{sec:proof-theo-prob}

\proof[Proof of \autoref{theo:main-prob}] 
Let $x_0 \in \R^d$, $r\in (0,\frac12)$, $x\in B_{r/4}(x_0)$. Using
\eqref{rem:potter} from \autoref{sec:appendix} with $\delta=1$, we see
that there is a constant
$c_0\geq 1$ so that 
\begin{align}\label{eq:potterp}
 \frac{L(s)}{L(s')}\leq c_0\left(\frac{s}{s'}\right)^{-\alpha-1},\quad
0<s<s'<1\,.
\end{align}

Define for $n\in \N$
\[
	r_n:=L^{-1}( L(\tfrac{r}{2})a^{n-1})\quad \text{ and }\quad s_n:=3 
\|u\|_\infty
b^{-(n-1)}
\]
for some constants $b\in (1,\frac{3}{2})$ and $a>c_0 2^{\alpha+1}$ that will be
chosen in the proof
independently of $n$, $r$ and $u$. As we explained in the
introduction, \autoref{theo:main-prob} trivially holds true of $\lim\limits_{r\to 0+}
L(r)$
is finite. Thus, we can assume $\lim\limits_{r\to 0+} L(r)$ to be infinite. This
implies that $r_n \to 0$ for $n \to \infty$ as it should be. 

We will use the following abbreviations:
\[
	B_n:=B_{r_n}(x), \quad \tau_n:=\tau_{B_n}, \quad m_n:=\inf_{B_n}u,
\quad M_n:=\sup_{B_n}u\,.
\]

We are going to prove 
 \begin{equation}\label{eq:thm-1}
 	M_k-m_k\leq s_k
 \end{equation}
 for all $k\geq 1$. 
 
Assume for a moment that \eqref{eq:thm-1} is proved. Then, for any $r\in
(0,\frac{1}{2})$ and $y\in B_{r/4}(x_0)\subset B_{r/2}(x)$ we can find $n\in
\N$ so that 
 \[
 	r_{n+1}\leq |y-x|< r_n\,.
 \]
Furthermore, since $L$ is decreasing, we obtain with 
$\gamma=\frac{\ln{b}}{\ln{a}}\in (0,1)$  
\begin{align*}
 	|u(y)-u(x)|&\leq s_n=3 b\|u\|_\infty a^{-n\frac{\ln b}{\ln a}} =3
b\|u\|_\infty \left[\frac{L(r_{n+1})}{
L(\frac{r}{2})}\right]^{-\frac{\ln{b}}{\ln{a}}} \leq 3
b\|u\|_\infty\left[\frac{L(|x-y|)}{L(\frac{r}{2})}\right]^{-\gamma}\,,
\end{align*}
which proves our assertion. Thus it remains to prove
(\ref{eq:thm-1}).

We are going to prove \eqref{eq:thm-1} by an inductive argument.
Obviously, $M_1-m_1\leq 2\|u\|_\infty\leq
s_1$. Since $1<b<\frac{3}{2}$, it follows that 
\[
M_2-m_2\leq 2\|u\|_\infty\leq 3\|u\|_\infty b^{-1}=s_2\,.
\]

Assume now that (\ref{eq:thm-1})  is true for all $k\in \{1,2,\ldots,n\}$ for 
some $n\geq 2$. 

Let $\varepsilon>0$ and take ${z_1},z_2 \in B_{n+1}$ so that
\[
	u({z_1})\leq m_{n+1}+\frac{\varepsilon}{2}\ \ \ \ \ \ \ \ u(z_2 )\geq
M_{n+1}-\frac{\varepsilon}{2}\,.
\]
It is enough to show that
\begin{equation}\label{eq:thm-2}
	u(z_2 )-u({z_1})\leq s_{n+1},
\end{equation}
since then 
\[
	M_{n+1}-m_{m+1}-\varepsilon\leq s_{n+1},
\]
which implies (\ref{eq:thm-1})  for $k=n+1$, since $\varepsilon>0$ was 
arbitrary. 

By the optional stopping theorem,
\begin{align*}
	u(z_2 )-u({z_1})=&\ \E_{z_2}[u(X_{\tau_n})-u({z_1})]\nonumber \\
	=&\ \E_{z_2}[u(X_{\tau_n})-u({z_1});X_{\tau_n}\in
B_{n-1}]\\&+\sum\limits_{i=1}^{n-2}\E_{z_2}[u(X_{\tau_n})-u({z_1});X_{\tau_n}\in
B_{n-i-1}\setminus B_{n-i}]\\&+\E_{z_2}[u(X_{\tau_n})-u({z_1});X_{\tau_n}\in
B_1^c]=I_1+I_2+I_3.
\end{align*}

Let $A=\{z \in B_{n-1}\setminus B_n|\, u(z)\leq \frac{m_n+M_n}{2}\}$. It is
sufficient to consider the case $\mu_{x}(A)\geq
\frac{1}{2}\mu_{x}(B_{n-1}\setminus B_n)$, where
$\mu_{x}$ is the measure defined by (\ref{eq:measure_mu}). In the
remaining case we would use 
$\mu_{x}((B_{n-1}\setminus B_n)\setminus A)\geq
\frac{1}{2}\mu_{x}(B_{n-1}\setminus B_n)$
 and could continue the proof with  $\|u\|_\infty -u$ and 
 \[
 	(B_{n-1}\setminus B_n)\setminus A=\left\{z\in B_{n-1}\setminus B_n
|\, \|u\|_\infty
-u(z)\leq
\frac{\|u\|_\infty
-m_n+\|u\|_\infty -M_n}{2}\right\}
 \]
instead of $u$ and $A$.

The estimate \eqref{eq:potterp} implies  
$a=\tfrac{L(r_{n+1})}{L(r_n)}\leq c_0(\tfrac{r_{n+1}}{r_n})^{-\alpha-1}$, from
where we deduce $r_{n+1}\leq r_n (c_0a^{-1})^{\frac{1}{\alpha+1}} \leq
\frac{r_n}{2}$ because of $a>c_0 2^{\alpha+1}$. Next, we make use of the
following property:
\begin{align}\label{eq:role_varphi} 
r_{n-1} = L^{-1}(L(\tfrac{r}{2}) a^{n-2}) = L^{-1}(\tfrac{1}{a}
L(\tfrac{r}{2}) a^{n-1}) =  L^{-1}(\tfrac{1}{a} L(r_n)) = \varphi_a(r_n) \,. 
\end{align}
Then by
\autoref{prop:hitting_general} (with
$r=r_n$ and $x_0=x$) we get  \[p_n:=\P_{z_2}(X_{\tau_n}\in A)\geq
C_5\frac{\ln{a}}{a}\,.\]

Hence,
\begin{align*}
 I_1&=\E_{z_2}[u(X_{\tau_n})-u({z_1});X_{\tau_n}\in B_{n-1}]\\
&=\E_{z_2}[u(X_{\tau_n})-u({z_1});X_{\tau_n}\in
A]+\E_{z_2}[u(X_{\tau_n})-u({z_1});X_{\tau_n}\in B_{n-1}\setminus
A]\\
&\leq \left(\tfrac{m_n+M_n}{2}-m_n\right)p_n+s_{n-1}(1-p_n)\\
&\leq \tfrac{1}{2}s_n p_n+s_{n-1}(1-p_n)\leq s_{n-1}(1-\tfrac{1}{2}p_n)\leq
s_{n-1}(1-\tfrac{C_5\ln{a}}{2a})\,.
\end{align*}

By \autoref{prop:prob-3}, 
\begin{align*}
	I_2& \leq \sum\limits_{i=1}^{n-2} s_{n-i-1}\P_{z_2}(X_{\tau_n}\not\in
B_{n-i})
	\leq C_4 \sum\limits_{i=1}^{n-2} s_{n-i-1}\tfrac{L(r_{n-i})}{L(r_n)}\\
	&\leq 3C_4\|u\|_\infty  \sum\limits_{i=1}^{n-2}
b^{-(n-i-2)}\tfrac{a^{n-i-1}}{a^{n-1}}\leq
3C_4 \|u\|_\infty \tfrac{b^{-n+3}}{a-b}\\
	&\leq C_4  \tfrac{b^3}{a-b}s_{n+1}\,.
\end{align*}
Similarly, by \autoref{prop:prob-3}, 
\[
	I_3\leq 2\|u\|_\infty\P_{z_2}(X_{\tau_n}\not\in B_1)\leq 2C_4
\|u\|_\infty
\tfrac{L(r_1)}{L(r_n)}=\tfrac{2C_4 }{3}
b\left(\tfrac{b}{a}\right)^{n-1}s_{n+1}\leq C_4 
\tfrac{b^2}{a} s_{n+1}\,.
\]
Hence,
\[
u(z_2)-u(z_1)\leq s_{n+1}b^2\left[1-\tfrac{C_5\ln{a}}{2a}+ \tfrac{C_4 b}{a-b}+
\tfrac{C_4 }{a}\right]\,.
\]

Since $a-b\geq \frac{a}{4}$ for $b\in (1,\frac{3}{2})$ and
$a>c_0 2^{\alpha+1}\geq 2$, it follows
that
\[q:=1-\tfrac{C_5\ln{a}}{2a}+
\tfrac{C_4  b}{a-b}+
\tfrac{C_4 }{a}\leq 1-\tfrac{C_5\ln{a}}{2a}+\tfrac{7C_4 }{a}=1-\tfrac{C_5\ln
a-14C_4 }{2a}\,.\]

Next, we choose $a>c_0 2^{\alpha+1}$ so large that $C_5\ln{a}-14C_4 >0$. Thus
$q<1$.
Finally, we choose $b\in
(1,\frac{3}{2})$ sufficiently small so that $b^2q<1$\,. 

Hence, (\ref{eq:thm-2}) holds, which finishes the proof of the inductive step 
and the theorem.
\qed

\begin{appendix}

\section{Slow and Regular Variation} \label{sec:appendix}

In this section we collect some properties of slowly resp. regularly varying
functions that are used in \autoref{sec:prob_estim} and 
\autoref{sec:proof-theo-prob}. 

\begin{definition}
 A measurable and positive function $\ell\colon (0,1)\rightarrow (0,+\infty)$ is
said to vary regularly at zero with index $\rho\in
\R$ if for every $\lambda>0$
\[
 \lim_{r\to 0+}\frac{\ell(\lambda r)}{\ell(r)}=\lambda^{\rho}\,.
\]
If a function varies regularly at zero with index $0$ it is said to vary
slowly at zero. For simplicity, we call such functions \emph{regularly
varying} resp. \emph{slowly varying} functions.
\end{definition}

Note that slowly resp. regularly varying functions include functions which are
neither
increasing nor decreasing. By \cite[Theorem 1.4.1 (iii)]{BGT87} it follows that
any function $\ell$ that
varies regularly with index $\rho\in \R$ is of the form
$\ell(r)=r^{\rho}\ell_0(r)$ for some
function $\ell_0$ that varies slowly. 

Assume $\int_0^1 s \, \ell(s) \, \d s \leq c$ for some $c >0$. Let $L\colon
(0,1)\rightarrow (0,+\infty)$ be defined by
\[
	L(r)=\int\limits_r^1\frac{\ell(s)}{s}\, \d s\ .
\]
The function $L$ is well defined because $L(r) = r^{-2} \int_r^1
r^2 \frac{\ell(s)}{s}\, \d s \leq r^{-2} \int_r^1
s \ell(s) \, \d s \leq c r^{-2}$. Note that \eqref{eq:K-Levy-bound} and
\eqref{eq:K3} imply that $\int_0^1 s \, \ell(s) \, \d s \leq c$ does hold in our
setting. We note that the function $L$ is always decreasing. 

Let us list further properties which are making use of in our proofs. Note that
they are established \cite{BGT87} for functions which are slowly resp. regularly
varying at the point $+\infty$. By a simple inversion we adopt the results
to functions which are slowly resp. regularly
varying at the point $0$.

\begin{enumerate}
 \item \label{rem:l_versus_L} If $\ell$ is slowly varying, then
\cite[Proposition
1.5.9a]{BGT87} $L$ is slowly varying  with 
\[\lim\limits_{r\to
0+}L(r)=+\infty\qquad \text{ and }\qquad \lim\limits_{r\to
0+}\frac{\ell(r)}{L(r)}=0\,.\]
\item \label{rem:karamata} If $\ell$ is slowly varying and $\rho>-1$, then
Karamata's theorem \cite[Proposition 1.5.8]{BGT87} ensures 
\[
  \lim_{r\to 0+}\frac{\int_0^r
s^{\rho}\ell(s)\d s}{r^{\rho+1}\ell(r)}=(\rho+1)^{-1}\,.
\]
\item \label{rem:l-reg_L-reg} If $\ell$ is
regularly varying of order $-\alpha < 0$ (in our case $0 < \alpha < 2$),
then
\cite[Theorem 1.5.11]{BGT87} 
\[
 \lim_{r\to 0+}\frac{L(r)}{\ell(r)}=\alpha^{-1}\,.
\] In particular, if $\ell$ is regularly varying of order $-\alpha < 0$,
then so is $L$. 
\item \label{rem:potter} Assume $\ell$ is regularly varying of order
$-\alpha \leq 0$ and stays bounded away from $0$ and $+\infty$ on every compact
subset of $(0,1)$. Then Potter's theorem \cite[Theorem 1.5.6 (ii)]{BGT87} 
implies
that for every $\delta>0$ there is a
constant
$C=C(\delta)\geq 1$ such that for $r,s \in (0,1)$ 
\begin{align*}
 \frac{\ell(r)}{\ell(s)}\leq C
\max\left\{ \left(\frac{r}{s}\right)^{-\alpha-\delta},
\left(\frac{r}{s}\right)^{-\alpha+\delta} \right\} \,.
\end{align*}
\item \label{rem:L_notzero} Since $L$ is non-increasing, we observe
$\lim\limits_{r\to 0+} L(r) \in (0,+\infty]$.
\end{enumerate}

\end{appendix}

%
\bibliographystyle{bibstyle_english_author-year2}
\bibliography{15-03-21_intrinsic-scaling} 

\end{document}